\documentclass[final,3p,numbers,sort]{elsarticle}
\usepackage{amssymb}
\usepackage{amsmath}
\usepackage{amsthm}
\usepackage{enumitem}
\usepackage{fancyhdr}
\usepackage{hyperref}
\hypersetup{%
  pdftitle={$L^p$ $(p\geq 1)$ solutions of general time interval multidimensional BSDEs},%
  pdfsubject={BSDEs with monotonicity generators},%
  pdfauthor={Lishun XIAO},%
  pdfkeywords={BSDEs, general time interval, monotonicity},%
  pdfstartview=FitH,
  CJKbookmarks=true,%
  bookmarksnumbered=true,%
  bookmarksopen=true,%
  colorlinks=true, linkcolor=blue, urlcolor=blue, citecolor=blue, %
}
\newcommand{\rtn}{\mathrm{\mathbf{R}}}
\newcommand{\N}{\mathrm{\mathbf{N}}}
\newcommand*{\PR}{\mathrm{\mathbf{P}}}
\newcommand*{\EX}{\mathrm{\mathbf{E}}}
\newcommand*{\dif}{\,\mathrm{d}}
\newcommand*{\F}{\mathcal{F}}
\newcommand{\s}{\mathcal{S}}
\newcommand{\M}{\mathrm{M}}
\newcommand*{\prs}{\dif\PR-\mathrm{a.s.}}
\newcommand*{\pts}{\dif\PR\times\dif t-\mathrm{a.e.}}
\newcommand{\Ito}{It\^{o}'s}
\newcommand{\tT}[1][0]{[#1,T]}
\newcommand{\intT}[2][T]{\int^{#1}_{#2}}
\newcommand{\me}{\mathrm{e}}
\newcommand{\one}[1]{{\bf 1}_{#1}}

\newtheorem{thm}{Theorem}
\newtheorem{lem}[thm]{Lemma}
\newtheorem{pro}[thm]{Proposition}
\newtheorem{rmk}[thm]{Remark}
\newtheorem{dfn}[thm]{Definition}
\newtheorem{ex}[thm]{Example}
\newproof{pf}{Proof}
\newproof{pou}{{\bf Proof of the uniqueness part}}
\newproof{poe}{{\bf Proof of the existence part}}

\makeatletter
\def\ps@pprintTitle{%
     \let\@oddhead\@empty
     \let\@evenhead\@empty
     \def\@oddfoot{}%
     \let\@evenfoot\@oddfoot}
\makeatother

\begin{document}
\begin{frontmatter}


\title{{\boldmath\bf $L^p$ $(p\geq 1)$ solutions of multidimensional BSDEs
with monotone generators in general time intervals}\tnoteref{found}}
\tnotetext[found]{Supported by the National Natural Science Foundation of China
(No. 11101422 and 11371362) and the Fundamental Research Funds for
the Central Universities (No. 2012LWB57).}
\date{April 23, 2012}
\author[cumt]{Lishun XIAO}
\author[cumt,fudan]{Shengjun FAN\corref{cor1}}%
\ead{f\_s\_j@126.com}
\cortext[cor1]{Corresponding author}

\author[cumt]{Na XU}

\address[cumt]{College of Science, China University of Mining and Technology,
  Xuzhou, Jiangsu, 221116, PR China}
\address[fudan]{School of Mathematical Sciences, Fudan University, Shanghai, 200433, PR China}

\begin{abstract}
In this paper, we are interested in solving general time interval
multidimensional backward stochastic differential equations in $L^p$ $(p\geq 1)$.
We first study the existence and uniqueness for $L^p$ $(p>1)$ solutions by the
method of convolution and weak convergence when the generator is monotonic in
$y$ and Lipschitz continuous in $z$ both non-uniformly with respect to $t$.
Then we obtain the existence and uniqueness for $L^1$ solutions with an additional
assumption that the generator has a sublinear growth in $z$ non-uniformly
with respect to $t$.
\end{abstract}

\end{frontmatter}


\section{Introduction}
\label{sec:Introduction}
In this paper, we consider the following multidimensional backward stochastic
differential equation (BSDE for short in the remaining):
\begin{equation}\label{eq:BSDEs}
  y_t=\xi+\intT{t} g(s,y_s,z_s)\dif s-\intT{t}z_s\dif B_s, \quad t\in\tT,
\end{equation}
where $T$ satisfies $0\leq T\leq +\infty$ called the terminal time;
$\xi$ is a $k$-dimensional random vector called the terminal condition;
the random function
$g(\omega,t,y,z):\Omega\times\tT\times\rtn^k\times\rtn^{k\times d}\mapsto\rtn^k$
is progressively measurable for each $(y,z)$, called the generator of BSDE
\eqref{eq:BSDEs}; and $B$ is a $d$-dimensional Brownian motion. The solution
$(y_t,z_t)_{t\in\tT}$ is a pair of adapted processes. The triple $(\xi,T,g)$ is
called the parameters of BSDE \eqref{eq:BSDEs}. We denote also by BSDE $(\xi,T,g)$
the BSDE with the parameters $(\xi,T,g)$.

The nonlinear case of multidimensional BSDEs has been introduced by
\citet{PardouxPeng1990SCL}. They proved an existence and uniqueness result
under the assumptions that the generator $g$ is uniformly Lipschitz continuous
in both $y$ and $z$. Their terminal time $T$ is a finite constant
and the terminal condition $\xi$ and the process
$\{g(t,0,0)\}_{t\in\tT}$ are square-integrable. Since then, the research of BSDEs
in both theory and application has been widely made by more and more people.
Many applications of BSDEs have been found in mathematical finance,
stochastic control, partial differential equations and so on (See \citet{ElKarouiPengQuenez1997MF}
for details). In both theory and application of BSDEs, it is essential to relax the Lipschitz
conditions on the generator $g$, improve the terminal time into the general
case and study the solutions under non-square integrable parameters.

Many works including \citet{Mao1995SPA,LepeltierSanMartin1997SPL,Kobylanski2000AP,Bahlali2001CRASPSI,Hamadene2003Bernoulli,BriandLepetierSanMrtin2007Bernoulli,Jia2008SPL,WangHuang2009SPL,FanJiangDavison2010CRASSI},
see also
the references therein, have weakened the Lipschitz condition on the generator $g$.
These works dealt only with the BSDEs with square-integrable parameters. But
the terminal condition $\xi$ and the process $\{g(t,0,0)\}_{t\in\tT}$
are not necessarily square-integrable in some practical applications.
Then $L^p$ $(p\geq1)$ solutions of BSDEs when $\xi$ and $\{g(t,0,0)\}_{t\in\tT}$
are $p$-integrable attracted a lot of attention of many researchers.
\citet{ElKarouiPengQuenez1997MF,BriandCarmona2000JAMSA,BriandDelyonHu2003SPA,Chen2010SAA,FanJiang2012JAMC},
for instance, proved some existence
and uniqueness results for $L^p$ $(p>1)$ solutions of BSDEs respectively in
different conditions. In particular, \citet{BriandDelyonHu2003SPA} proved
the existence and uniqueness for $L^p$ $(p>1)$ solutions of multidimensional BSDEs
when the generator $g$ is monotonic in $y$ and is Lipschitz continuous in $z$;
and with an additional assumption that $g$ has a kind of sublinear growth in $z$
they obtained the existence and uniqueness for $L^1$ solutions.
To our knowledge, only few papers solved BSDEs with only integrable parameters
besides \citet{BriandDelyonHu2003SPA},
such as \citet{Peng1997BSDEP} and \citet{FanLiu2010SPL}. Particularly,
\citet{FanLiu2010SPL} obtained the existence and uniqueness for $L^1$ solutions
of one-dimensional BSDEs when the generator $g$ is Lipschitz continuous
in $y$ and $\alpha$-H\"older $(0<\alpha<1)$ continuous in $z$, which may be a
basic result for $L^1$ solutions of BSDEs. However, all these works talked above
dealt only with the BSDEs with a finite time interval. BSDEs with general
time intervals have not been researched widely.

In the sequl, we introduce some papers which studied the existence and uniqueness for
solutions of BSDEs with general time intervals.
\citet{ChenWang2000JAMSA} first improved the terminal time into the general
case and proved
the existence and uniqueness for $L^2$ solutions of one-dimensional BSDEs by the
fixed point theorem under the assumptions that the generator $g$ is Lipschitz
continuous in $(y,z)$ non-uniformly with respect to $t$, which actually extended
the result of \citet{PardouxPeng1990SCL} into the general time interval case.
\citet{FanJiang2010SPL} and \citet{FanJiangTian2011SPA} respectively relaxed the Lipschitz condition
of \citet{ChenWang2000JAMSA}, and obtained the existence and uniqueness
result for $L^2$ solutions of BSDEs with general time intervals.
Recently, \citet{FanJiang2011Stochastics} investigated the existence and uniqueness
for $L^p$ $(p>1)$ solutions of multidimensional BSDEs with general
time intervals under some weaker assumptions. However, all these works
need a linear-growth condition of the generator $g$ with respect to $y$ to guarantee
the existence of $L^p$ $(p>1)$ solutions. On the other hand, to our knowledge,
there are no papers which have studied the $L^1$ solutions of BSDEs with general
time intervals.

In this paper, under a monotonicity condition and a general growth condition for
the generator $g$ with respect to $y$ we establish a general existence and uniqueness
result for $L^p$ $(p\geq 1)$ solutions of multidimensional BSDEs with
general time intervals (see Theorem \ref{thm:MainResultLpSolution}
in Section \ref{sec:LpSolution} and Theorem
\ref{thm:MainResultL1Solution} in Section \ref{sec:IntegrableParameters}).
In particular, the first part of this paper is devoted to proving the existence
and uniqueness for $L^p$ $(p>1)$ solutions when the generator $g$ is monotonic
and has a general growth in $y$ and is Lipschitz continuous in $z$, which are
both non-uniform with respect to $t$ (see \ref{H:gGeneralizedGeneralGrowthInY} --
\ref{H:gLipschitzInZ} in Section \ref{sec:LpSolution}).
After that, we study the existence and uniqueness for $L^1$ solutions under the same conditions together
with an additional sublinear growth assumption in $z$ (see \ref{H:gSublinearGrowthInZ} in
Section \ref{sec:IntegrableParameters}). Note that the $u(t)$, $v(t)$ and $\gamma(t)$
appearing in assumptions \ref{H:gMonotonicityInY} -- \ref{H:gSublinearGrowthInZ}
may be unbounded and their integrability is the only requirement (see Remarks
\ref{rmk:UVUnbounded} and \ref{rmk:ComparisonWithBraindResults} for details). Our results
actually extend and improve the results of \citet{BriandDelyonHu2003SPA} into the
general time interval case when the assumptions on the generator $g$ is not
necessarily uniform with respect to $t$. Besides, our results also include the
corresponding results of \citet{PardouxPeng1990SCL}, \citet{Pardoux1999NADEC} and
\citet{ChenWang2000JAMSA} as its particular cases.

The rest of this paper is organized as follows. Section \ref{sec:PreliminariesAndPriorEstimates}
introduces some notations and lemmas used in the whole paper, and also establishes
some important apriori estimates for solutions of BSDE \eqref{eq:BSDEs}.
Section \ref{sec:LpSolution} puts forward and proves the
existence and uniqueness result for the $L^p$ $(p>1)$ solutions.
Section \ref{sec:IntegrableParameters} shows the
existence and uniqueness for the $L^1$ solutions. \ref{sec:Appendix} gives some detailed
proofs of lemmas.

\section{Preliminaries and apriori estimates}
\label{sec:PreliminariesAndPriorEstimates}
Although many researchers use the same notations in studying BSDEs, we will
still introduce the following notations in order to make the paper easy to read.

First of all, let $(\Omega,\F,\PR)$ be a probability space carrying a standard
$d$-dimensional Brownian motion $(B_t)_{t\geq 0}$ and let $(\F_t)_{t\geq 0}$
be the natural $\sigma$-algebra filtration generated by $(B_t)_{t\geq 0}$.
We assume that $\F_T=\F$ and $(\F_t)_{t\geq 0}$ is right-continuous and complete.
In this paper, the Euclidean norm of a vector $y\in \rtn^k$ will be defined by
$|y|$, and for a $k\times d$ matrix $z$, we define $|z|=\sqrt{Tr(zz^*)}$,
where $z^*$ is the transpose of $z$. Let $\langle x,y\rangle$ represent the
inner product of $x$, $y\in\rtn^k$.

For each real number $p>0$, let $L^p(\Omega,\F_T,\PR;\rtn^k)$ be the set of
$\rtn^k$-valued and $\F_T$-measurable random variables $\xi$ such that
$\|\xi\|^p_{L^p}:=\EX[|\xi|^p]<+\infty$ and let ${\s}^p(0,T;\rtn^k)$ (or $\s^p$ for notation
convenience) denote the set of $\rtn^k$-valued, adapted and continuous processes
$(Y_t)_{t\in\tT}$ such that
$$\|Y\|_{{\s}^p}:=
  \left(\EX
    \left[
       \sup_{t\in\tT}|Y_t|^p
    \right]
  \right)^{1\wedge 1/p}<+\infty. $$
If $p\geq 1$, $\|\cdot\|_{{\s}^p}$ is a norm on ${\s}^p$ and if $p\in (0,1)$,
$(Y,Y')\mapsto \|Y-Y'\|_{{\s}^p}$ defines a distance on ${\s}^p$. Under
this metric, ${\s}^p$ is complete. Moreover, let $\M^p(0,T;\rtn^{k\times d})$
(or $\M^p$ for notation convenience) denote the set of
(equivalent classes of) $(\F_t)$-progressively measurable
${\rtn}^{k\times d}$-valued processes $(Z_t)_{t\in\tT}$ such that
$$\|Z\|_{{\M}^p}:=
  \left\{\EX
    \left[
      \left(
        \int_0^T |Z_t|^2\dif t
      \right)^{p\over 2}
    \right]
  \right\}^{1\wedge 1/p}<+\infty.$$
For any $p\geq 1$, $\M^p$ is a Banach space endowed with this norm and for any
$p\in (0,1)$, $\M^p$ is a complete metric space with the resulting distance.
We also denote by $\|\cdot\|_{\s^p\times\M^p}$ the norm in the space $\s^p\times \M^p$
for any $p>1$ with the following definition
  \begin{equation*}
    \|(y,z)\|_{\s^p\times\M^p}:=
    \left(\EX
      \left[
        \sup_{t\in\tT}|y_t|^p+\left(\intT{0}|z_t|^2\dif t\right)^{p\over 2}
      \right]
    \right)^{1\over p}.
  \end{equation*}

Let us recall that a continuous process $(Y_t)_{t\in\tT}$ belongs to the
class (D) if the family
$\{Y_\tau:\tau\in\Sigma_T\}$
is uniformly integrable, where $\Sigma_T$ stands for the set of
all stopping times $\tau$ such that $\tau\leq T$.
For a process $(Y_t)_{t\in\tT}$ belonging to the class (D), we define
$$\|Y\|_1=\sup\{\EX[|Y_\tau|]:\tau\in\Sigma_T\}.$$
The space of $(\F_t)$-progressively measurable continuous
processes which belong to the class (D) is complete under this norm.

As mentioned above, we will deal only with the multidimensional BSDE which is
an equation of type \eqref{eq:BSDEs}, where the
terminal condition $\xi$ is $\F_T$-measurable, the terminal time $T$ satisfies
$0\leq T\leq +\infty$ and the generator $g$ is $(\F_t)$-progressively
measurable for each $(y, z)$.

\begin{dfn}
  Let $T$ satisfy $0\leq T\leq+\infty$.
  A pair of processes $(y_t,z_t)_{t\in\tT}$ taking values in $\rtn^k\times\rtn^{k\times d}$
  is called a solution of BSDE \eqref{eq:BSDEs},
  if $(y_t,z_t)_{t\in\tT}$ is $(\F_t)$-adapted
  and satisfies that $\prs$, $t\mapsto y_t$ is continuous, $t\mapsto z_t$ belongs
  to $L^2(0,T)$, $t\mapsto g(t,y_t,z_t)$ belongs to $L^1(0,T)$ and $\prs$,
  BSDE \eqref{eq:BSDEs} holds true for each $t\in\tT$.
\end{dfn}

Let us introduce the following Lemma \ref{lem:chenwang} which
comes from Lemma 1.1 and Theorem 1.2 of \citet{ChenWang2000JAMSA}. Note that
Lemma 2
holds also in the multidimensional case since the proofs of Lemma 1.1 and Theorem
1.2 of Chen and Wang are done via a standard contraction argument
combined with apriori estimates without using comparison theorem.

\begin{lem}\label{lem:chenwang}
Assume that $0\leq T\leq +\infty$, $\xi\in L^2(\Omega,\F_T,\PR;\rtn^k)$ and
the following assumptions hold:
\begin{enumerate}
  \item[(C1)] $\EX\left[\left(\intT{0}|g(t,0,0)|\dif t\right)^2\right]<+\infty$;
  \item[(C2)] There exist two deterministic functions
              $u(t), v(t):[0,T]\mapsto\rtn^+$ with $\intT{0}\big(u(t)+v^2(t)\big)\dif t<+\infty$
              such that $\pts$, for each $(y_i,z_i)\in\rtn^k\times\rtn^{k\times d}$,
              $i=1,2$,
              $$|g(t,y_1,z_1)-g(t,y_2,z_2)|\leq u(t)|y_1-y_2|+v(t)|z_1-z_2|.$$
\end{enumerate}
Then BSDE \eqref{eq:BSDEs} has a unique solution in the space $\s^2\times \M^2$.
\end{lem}

Throughout this paper we will use the Corollary 2.3 in \citet{BriandDelyonHu2003SPA}
several times. So we list it as a lemma. Note that this conclusion holds still
true for $T=+\infty$.
\begin{lem}\label{lem:Tanaka-Briand}
  If $(y_t,z_t)_{t\in\tT}$ is a solution of BSDE \eqref{eq:BSDEs}, $p\geq 1$,
  $c(p)=p[(p-1)\wedge 1]/2$ and $0\leq t\leq u\leq T\leq +\infty$, then
  \begin{align*}
    |y_t|^p+c(p)\intT[u]{t}|y_s|^{p-2}\one{|y_s|\neq 0}|z_s|^2\dif s
    \leq&\ |y_u|^p+p\intT[u]{t}|y_s|^{p-2}\one{|y_s|\neq 0}
      \langle y_s,g(s,y_s,z_s)\rangle\dif s\\
    &-p\intT[u]{t}|y_s|^{p-2}\one{|y_s|\neq 0}\langle y_s,z_s\dif B_s\rangle.
  \end{align*}
\end{lem}

Next we will establish some apriori estimates which play an important role
in proving our main results. In stating them, it is useful to introduce the
following assumption on the generator $g$, where $p>0$ and $0\leq T\leq+\infty$.
\begin{enumerate}
  \item[(A)] There exist two nonnegative functions $\mu(t)$,
        $\lambda(t):\tT\mapsto\rtn^+$ with $\intT{0}\big(\mu(t)+\lambda^2(t)\big)\dif t<+\infty$ such that
        $\pts$, for each $(y,z)\in\rtn^k\times\rtn^{k\times d}$,
        \[
          \langle y,g(t,y,z)\rangle\leq \mu(t)|y|^2+\lambda(t)|y||z|+f_t|y|,
        \]
        where $(f_t)_{t\in\tT}$ is a nonnegative and $(\F_t)$-progressively
        measurable processes with $\EX[(\intT{0}f_t\dif t)^p]<+\infty$.
\end{enumerate}

The following Lemmas \ref{lem:estimate-z} and \ref{lem:estimate-y}
give some estimates for $L^p$ solutions of BSDE
\eqref{eq:BSDEs} with $0\leq T\leq +\infty$ in the spirit of the
work in \citet{BriandDelyonHu2003SPA}, whose proofs are given
in \ref{sec:Appendix}.

\begin{lem}\label{lem:estimate-z}
  Assume that $0\leq T\leq+\infty$, $g$ satisfies assumption (A),
  $(y_t,z_t)_{t\in\tT}$ is a solution of BSDE \eqref{eq:BSDEs} such that
  $(y_t)_{t\in\tT}\in\s^p$ with $p>0$, $\beta(t):\tT\mapsto\rtn^+$ with
  $\intT{0}\beta(t)\dif t<+\infty$ and $\beta(t)\geq 2\big(\mu(t)+\lambda^2(t)\big)$.
  Then there exists a constant $C^1_p>0$ depending only on
  $p$ such that for each $0\leq r\leq t\leq T$,
  \begin{align*}
    \EX\!\left[\!\left.
       \left(
         \intT{t}\me^{\intT[s]{0}\beta(\bar u)\dif \bar u}|z_s|^2\dif s
       \right)^{p\over 2}\right|\F_r
       \right]\!
    \leq C^1_p\EX\!
       \left[\left.
         \sup_{s\in\tT[t]}
         \left(
           \me^{\frac{p}{2}\intT[s]{0}\beta(\bar u)\dif \bar u}|y_s|^p
         \right)\!
       +\!\left(\intT{t}\me^{\frac{1}{2}\intT[s]{0}\beta(\bar u)\dif \bar u}f_s\dif s\right)^p
       \right|\F_r
       \right].
  \end{align*}
\end{lem}

\begin{lem}\label{lem:estimate-y}
  Let the assumptions of Lemma 4 hold and assume further that $p>1$ and
  $\beta(t)\geq 2\{\mu(t)+\lambda^2(t)/[1\wedge(p-1)]\}$. Then there exists
  a constant $C^2_p>0$ depending only on $p$ such that for each $0\leq r\leq t\leq T$,
  \begin{align*}
   \EX\left[\left.
       \sup_{s\in\tT[t]}
       \left(
         \me^{\frac{p}{2}\intT[s]{0}\beta(\bar u)\dif \bar u}|y_s|^p
       \right)\right|\F_r
     \right]
   \leq C^2_p\EX
     \left[
       \left.\me^{\frac{p}{2}\intT{0}\beta(\bar u)\dif \bar u}|\xi|^p
       +\left(\intT{t}\me^{\frac{1}{2}\intT[s]{0}\beta(\bar u)\dif \bar u}f_s\dif s\right)^p\right|\F_r
     \right].
  \end{align*}
\end{lem}
\vspace{0.5em}

Combing Lemma \ref{lem:estimate-z} and Lemma \ref{lem:estimate-y} we can obtain the
following Proposition \ref{pro:estimate}.

\begin{pro}\label{pro:estimate}
  Let the assumptions in Lemma \ref{lem:estimate-y} hold and $p>1$, then there
  exists a constant $C_p>0$ depending only on $p$ such that for each
  $0\leq r\leq t\leq T$,
  \begin{align*}
   & \EX\!\!
     \left[\!\left.
        \intT{t}\!\!\me^{\frac{p}{2}\intT[s]{0}\beta(\bar u)\dif \bar u}
        \beta(s)|y_s|^p\dif s\right|\F_r
     \right]\!+\!
     \EX\!\!
     \left[\!\left.
       \sup_{s\in\tT[t]}\!
       \left(\me^{\frac{p}{2}\intT[s]{0}\beta(\bar u)\dif \bar u}|y_s|^p
       \right)\right|\F_r
     \right]\!+\!
     \EX\!\!
     \left[\!\left.
      \left(
        \intT{t}\!\!\me^{\intT[s]{0}\beta(\bar u)\dif \bar u}|z_s|^2\dif s
      \right)^{p\over 2}\right|\F_r
     \right]\\
   & \leq C_p
    \left\{
      \EX
      \left[\left.
        \me^{\frac{p}{2}\intT{0}\beta(\bar u)\dif \bar u}|\xi|^p\right|\F_r
      \right]+\EX
      \left[\left.
        \left(
          \intT{t}\me^{\frac{1}{2}\intT[s]{0}\beta(\bar u)\dif \bar u}f_s\dif s
        \right)^p\right|\F_r
      \right]
    \right\}.
  \end{align*}

\end{pro}

\section{$L^p$ $(p>1)$ solution}
\label{sec:LpSolution}
This section will give an existence and uniqueness result for
$L^p$ $(p>1)$ solutions of BSDE \eqref{eq:BSDEs} with $0\leq T\leq +\infty$
under the assumptions that the generator $g$ is monotonic and has a
general growth in $y$, and is Lipschitz continuous in $z$, which are both
non-uniform with respect to $t$.

First, we introduce the following assumptions with respect
to the generator $g$ of BSDE \eqref{eq:BSDEs} where $p>1$ and $0\leq T\leq +\infty$.
In stating them we always suppose that $u(t)$, $v(t):\tT\mapsto \rtn^+$
are two deterministic functions such that $\intT{0}\big(u(t)+v^2(t)\big)\dif t<+\infty$.

\begin{enumerate}
\renewcommand{\theenumi}{(H\arabic{enumi})}
\renewcommand{\labelenumi}{\theenumi}

  \item\label{H:p-integrable}
          $\EX
          \left[
            \left(
              \intT{0}|g(t,0,0)|\dif t
            \right)^p
          \right]<+\infty$;
  \item\label{H:gContinuousInY}
        $\pts$, for each $z\in\rtn^{k\times d}$, $y\mapsto g(t,y,z)$
        is continuous;

  \item\label{H:gGeneralizedGeneralGrowthInY}
        $g$ has a general growth in $y$, i.e., for each $r'\in\rtn^+$,
        we have
        $$
          \psi_{r'}(t):=\sup_{|y|\leq r'}|g(t,y,0)-g(t,0,0)|\in L^1(\tT\times\Omega);
        $$
  \item\label{H:gMonotonicityInY} $g$ is monotonic in $y$ non-uniformly with
        respect to $t$, i.e.,
        $\pts$, for each $y_1,y_2\in\rtn^k$, $z\in\rtn^{k\times d}$, we have
        \[
          \langle y_1-y_2,g(t,y_1,z)-g(t,y_2,z)\rangle\leq u(t)|y_1-y_2|^2;
        \]

  \item\label{H:gLipschitzInZ} $g$ is Lipschitz continuous in $z$ non-uniformly
        with respect to $t$, i.e., $\pts$,
        for each $y\in\rtn^k$, $z_1,z_2\in\rtn^{k\times d}$, we have
        \[
          |g(t,y,z_1)-g(t,y,z_2)|\leq v(t)|z_1-z_2|.
        \]
\end{enumerate}
Moreover, we need the following assumption.
\begin{enumerate}
  \renewcommand{\theenumi}{(H\arabic{enumi}')}
  \renewcommand{\labelenumi}{\theenumi}
  \setcounter{enumi}{2}
  \item\label{H:gGeneralGrowthInY} There
        exists a continuous increasing function $\varphi:\rtn^+\mapsto\rtn^+$
        such that $\pts$, for each $y\in\rtn^k$, $z\in\rtn^{k\times d}$, we have
       \[
         |g(t,y,z)|\leq |g(t,0,z)|+u(t)\varphi(|y|).
       \]
\end{enumerate}

\begin{rmk}\label{rmk:UVUnbounded}
  It is clear that assumption \ref{H:gGeneralizedGeneralGrowthInY} is
  weaker than assumption \ref{H:gGeneralGrowthInY}.
  In addition, it should be noted that in the corresponding assumptions
  of \citet{BriandDelyonHu2003SPA} and \citet{Pardoux1999NADEC}, the $u(t)$ and
  $v(t)$ in \ref{H:gGeneralGrowthInY}, \ref{H:gMonotonicityInY} and
  \ref{H:gLipschitzInZ} are assumed to be bounded by a constant $c>0$
  since they work with continuous functions in a compact time interval.
  In our framework, they may be unbounded.
\end{rmk}

\begin{rmk}\label{rmk:ChangeOfVariables}
  It is not hard to verify that $(y_t,z_t)_{t\in\tT}$ is a solution of BSDE $(\xi,T,g)$ iff
  \[
    (\overline{y}_t,\overline{z}_t):=
    \left(\me^{\intT[t]{0}u(s)\dif s}y_t,\me^{\intT[t]{0}u(s)\dif s}z_t\right)
  \]
  is a solution of BSDE $(\me^{\intT{0}u(s)\dif s}\xi,T,\overline{g})$, where
  \[
    \overline{g}(t,y,z):=\me^{\intT[t]{0}u(s)\dif s}g\big(t,\me^{-\intT[t]{0}u(s)\dif s}y,
    \me^{-\intT[t]{0}u(s)\dif s}z\big)-u(t)y.
  \]
  We can check that $\overline g$ satisfies the previous assumptions as $g$,
  but with \ref{H:gMonotonicityInY} replaced by
  \begin{enumerate}[itemindent=0.5ex]
  \renewcommand{\theenumi}{(H\arabic{enumi}')}
  \renewcommand{\labelenumi}{\theenumi}
  \setcounter{enumi}{3}
    \item\label{H:gMonotonicityInYWithU=0}
    $\langle y_1-y_2,g(t,y_1,z)-g(t,y_2,z)\rangle\leq 0.$
  \end{enumerate}
  Therefore, without loss of generality, we can assume that $g$ satisfies
  \ref{H:gMonotonicityInYWithU=0} provided that $g$ satisfies
  \ref{H:gMonotonicityInY}.
\end{rmk}

The main result of this section is as follows.

\begin{thm}\label{thm:MainResultLpSolution}
  Assume that $0\leq T\leq +\infty$, $p>1$ and $g$ satisfies assumptions
  \ref{H:p-integrable} -- \ref{H:gLipschitzInZ}. Then for each
  $\xi\in L^p(\Omega,\F_T,\PR;\rtn^k)$, BSDE \eqref{eq:BSDEs}
  has a unique solution $(y_t,z_t)_{t\in\tT}$ in
  $\s^p\times \M^p$.
\end{thm}

Next we will prove Theorem \ref{thm:MainResultLpSolution}.
By Remark \ref{rmk:ChangeOfVariables} we shall always assume that $g$
satisfies \ref{H:p-integrable} --
\ref{H:gGeneralizedGeneralGrowthInY}, \ref{H:gMonotonicityInYWithU=0}
and \ref{H:gLipschitzInZ}. Let us prove the uniqueness part
first and then the existence part.
\begin{proof}[\bf Proof of the uniqueness part of Theorem \ref{thm:MainResultLpSolution}]
  Let $(y_t^1,z_t^1)_{t\in\tT}$ and $(y_t^2,z_t^2)_{t\in\tT}$ be two solutions
  of BSDE \eqref{eq:BSDEs} such that both $(y^1_t,z^1_t)_{t\in\tT}$ and
  $(y^2_t,z^2_t)_{t\in\tT}$ belong to $\s^p\times\M^p$. We set
  $\hat{y}:=y^1-y^2$ and $\hat{z}:=z^1-z^2$, then
  $(\hat{y}_t,\hat{z}_t)_{t\in\tT}$ is a solution of the
  following BSDE in $\s^p\times\M^p$,
  \begin{equation}\label{eq:BSDE-hatg}
    \hat{y}_t=\intT{t}\hat{g}(s,\hat{y}_s,\hat{z}_s)\dif s-\intT{t}\hat{z}_s\dif B_s,\quad t\in\tT,
  \end{equation}
  where $\hat{g}(t,y,z):=g(t,y+y^2_t,z+z^2_t)-g(t,y^2_t,z^2_t)$
  for each $(y,z)\in\rtn^k\times\rtn^{k\times d}$. It follows from
  assumptions \ref{H:gMonotonicityInYWithU=0} and \ref{H:gLipschitzInZ} on $g$ that
  \begin{align*}
    \langle y,\hat{g}(t,y,z)\rangle
    &=\langle y,g(t,y+y^2_t,z+z^2_t)-g(t,y^2_t,z+z^2_t)\rangle
     +\langle y,g(t,y^2_t,z+z^2_t)-g(t,y^2_t,z^2_t)\rangle\\
    &\leq \langle y, g(t,y^2_t,z+z^2_t)-g(t,y^2_t,z^2_t)\rangle
      \leq v(t)|y||z|,
  \end{align*}
  which means that assumption (A) is satisfied for the generator
  $\hat{g}(t,y,z)$ of BSDE \eqref{eq:BSDE-hatg} with $\mu(t)\equiv 0$,
  $\lambda(t)=v(t)$ and $f_t\equiv 0$. Thus, by Proposition \ref{pro:estimate}
  with $r=t=0$ we know that
  \begin{equation*}
    \EX\left[
         \sup_{s\in\tT}|\hat{y}_s|^p
         +\left(\intT{0}|\hat{z}_s|^2\dif s\right)^{p\over 2}
       \right]=0.
  \end{equation*}
  Therefore, $(\hat{y}_t,\hat{z}_t)_{t\in\tT}=(0,0)$. The
  uniqueness part is then complete.
\end{proof}

Next we begin to prove the existence part of Theorem \ref{thm:MainResultLpSolution}.
The proof method is enlightened by \citet{Pardoux1999NADEC}, \citet{BriandDelyonHu2003SPA},
and \citet{BriandCarmona2000JAMSA}.
More precisely, the techniques
applied in the first step, the convolution and the weak convergence, are lent from
\citet{Pardoux1999NADEC}, and the truncation techniques are taken partly from
\citet{BriandDelyonHu2003SPA} and \citet{BriandCarmona2000JAMSA}. It should be mentioned that
since we have changed the terminal time from the finite case to the general case,
the space of $(y_t)_{t\in\tT}$ from the square-integrable space to $\s^p(0,T;\rtn^k)$,
and the $p$-integrable condition of $\{g(t,0,0)\}_{t\in\tT}$ from
$\EX[\intT{0}|g(t,0,0)|^p\dif t]<+\infty$ to
$\EX[(\intT{0}|g(t,0,0)|\dif t)^p]<+\infty$, some new troubles come up naturally
when we combine those techniques mentioned above together. For example,
in the case of $T=+\infty$, the integration of a constant over $\tT$ is not
finite anymore; $\|X\|_{\M^p}\leq C\|X\|_{\s^p}$ may not hold any longer;
and the condition $\intT{0}v^2(s)\dif s<+\infty$ can not imply $\intT{0}v(s)\dif s<+\infty$.
Additionally, from the point of technique view, in order to prove the existence part of
Theorem \ref{thm:MainResultLpSolution}, we need an existence result under
assumption \ref{H:gGeneralGrowthInY}, but it has not been proved
in the general time interval case. All these troubles will be solved using
different procedures.

\begin{proof}[\bf Proof of the existence part of Theorem \ref{thm:MainResultLpSolution}]
The proof will be done by four steps as follows:
\begin{itemize}
  \item With the help of Lemma \ref{lem:chenwang} and Proposition \ref{pro:estimate},
        by applying an approximation method via convolution smoothing as well as an argument
        on weak convergence borrowed from \citet{Pardoux1999NADEC}, we prove
        the existence of a solution in $\s^2\times\M^2$ for the following BSDE:
        \begin{equation}\label{eq:BSDE-g-yt-Vt}
          y_t=\xi+\intT{t}g(s,y_s,V_s)\dif s-\intT{t}z_s\dif B_s,\quad t\in\tT,
        \end{equation}
        under assumptions \ref{H:gContinuousInY}, \ref{H:gGeneralGrowthInY},
        \ref{H:gMonotonicityInYWithU=0} and \ref{H:gLipschitzInZ}, provided that
        $V\in\M^p$ and there exists a nonnegative constant
        $K$ such that
        \begin{equation}\label{eq:AssumptionFirstStep}
          |\xi|\leq K,\prs \quad \text{and} \quad |g(t,0,V_t)|\leq K\me^{-t},\pts.
        \end{equation}
  \item By using a particular truncation technique, we prove that the assumption
        \ref{H:gGeneralGrowthInY} in the above step can be weakened to \ref{H:gGeneralizedGeneralGrowthInY}.
  \item With the help of Proposition \ref{pro:estimate}, by a similar truncation
        argument to that in \citet{BriandDelyonHu2003SPA}, we prove that for
        each $\xi\in L^p(\Omega,\F_T,\PR;\rtn^k)$ and $V\in\M^p$,
        BSDE \eqref{eq:BSDE-g-yt-Vt} has a solution in $\s^p\times\M^p$ under
        assumptions \ref{H:p-integrable} -- \ref{H:gGeneralizedGeneralGrowthInY},
        \ref{H:gMonotonicityInYWithU=0} and \ref{H:gLipschitzInZ}.
  \item We construct a strict contraction by subdividing the time interval $\tT$ to
        show the existence of a solution to BSDE \eqref{eq:BSDEs} in the space
        $\s^p\times\M^p$ for each $\xi\in L^p(\Omega,\F_T,\PR;\rtn^k)$ under
        assumptions \ref{H:p-integrable} -- \ref{H:gGeneralizedGeneralGrowthInY},
        \ref{H:gMonotonicityInYWithU=0} and \ref{H:gLipschitzInZ}, which is the
        desired result.
\end{itemize}

On the whole, the first three steps deal with the case where the generator $g$
is independent of $z$, and the last step considers the general case.

{\medskip\noindent\bf First step: } 
Now we assume that $\xi\in L^p(\Omega,\F_T,\PR;\rtn^k)$, $V\in \M^p(0,T;\rtn^{k\times d})$
and that \ref{H:gContinuousInY}, \ref{H:gGeneralGrowthInY}, \ref{H:gMonotonicityInYWithU=0},
\ref{H:gLipschitzInZ} and \eqref{eq:AssumptionFirstStep} hold true.
For notational convenience, in this step we set, for each $y\in\rtn^k$,
\[f(t,y):=g(t,y,V_t).\]
Clearly, we have that
\[
  \EX
  \left[
    |\xi|^2
  \right]<+\infty \quad \text{and}\quad
  \EX
  \left[
    \left(
      \intT{0}|f(t,0)|\dif t
    \right)^2
  \right]<+\infty.
\]

Let $\rho_n(x):=n^k\rho(nx)$, where $\rho: \rtn^k\mapsto\rtn^+$ is a nonnegative
$\mathcal{C}^\infty$ function with the unit ball for compact support and which
satisfies $\int_{\rtn^k}\rho(x)\dif x=1$. We define for each
$(\omega,t,y)\in\Omega\times\tT\times\rtn^k$,
\begin{equation}\label{eq:g-n-convolution}
  f_n(t,y):=\big(\rho_n(\cdot)*f(t,\cdot)\big)(y)
  =\int_{\rtn^k}\rho_n(x)f(t,y-x)\dif x.
\end{equation}
Then $f_n$ is a $(\F_t)$-progressively measurable process for each $y\in\rtn^k$ and
\begin{equation}\label{eq:g-n-convolution-equivalent}
  f_n(t,y)=\int_{\rtn^k}\rho(x)f(t,y-\frac{x}{n})\dif x
  =\int_{\{x:|x|\leq 1\}}\rho(x)f(t,y-\frac{x}{n})\dif x.
\end{equation}
Concerning $f_n(t,y)$, we have the following Lemma \ref{lem:AppendixProofOfgnWithLocallyLip}, whose
proof is given in \ref{sec:Appendix}.
\begin{lem}\label{lem:AppendixProofOfgnWithLocallyLip}
  Take driver $g$ under \ref{H:gGeneralGrowthInY} -- \ref{H:gMonotonicityInYWithU=0}
  and define $\phi(u'):=\varphi(u'+1)$ for each $u'\in\rtn^+$. Then for each $n\in\N$,
  $f_n(t,y)$ satisfies \ref{H:gContinuousInY},
  \ref{H:gGeneralGrowthInY} with $\varphi$ replaced by $\phi$, and
  \ref{H:gMonotonicityInYWithU=0}. Furthermore, we have
  \begin{equation}\label{eq:FnEnlargedByEAndU}
    |f_n(t,0)|\leq K\me^{-t}+u(t)\phi(0),
  \end{equation}
  and for each $y_1$, $y_2\in\rtn^k$ with $|y_1|\leq m$ and $|y_2|\leq m$, there
  exists a constant $C^n_m$ depending on $m$ and $n$ such that
  \begin{equation}\label{eq:fnLocallyLipschitzContinuousInY}
    |f_n(t,y_1)-f_n(t,y_2)|\leq C^n_m\big(\me^{-t}+u(t)\big)|y_1-y_2|,
  \end{equation}
  which means that $f_n(t,y)$ is locally Lipschitz continuous in
  $y$ non-uniformly with respect to $t$.
\end{lem}

  We now define for each $q\in \N$,
  \begin{equation*}
    f_{n,q}(t,y):=f_n(t,\pi_q(y)),
  \end{equation*}
  where and hereafter for each $u'\in\rtn^+$ and $x\in\rtn^k$,
  \[
    \pi_{u'}(x):=\frac{u' x}{u'\vee |x|}.
  \]
  By \eqref{eq:FnEnlargedByEAndU} we know that
  $f_{n,q}(t,0)=f_n(t,0)$ satisfies (C1) in Lemma \ref{lem:chenwang}.
  Furthermore, it follows from \eqref{eq:fnLocallyLipschitzContinuousInY} that
  there exists a constant $K_{n,q}$ depending on $n$ and $q$ such that
  for each $y_1$, $y_2\in\rtn^k$, in view of $|\pi_q(y)|\leq q$ for each $y\in\rtn^k$
  and $|\pi_q(y_1)-\pi_q(y_2)|\leq |y_1-y_2|$,
  \begin{equation*}
    |f_{n,q}(t,y_1)-f_{n,q}(t,y_2)|\leq K_{n,q}\big(\me^{-t}+u(t)\big)|y_1-y_2|,
  \end{equation*}
  which implies that $f_{n,q}$ satisfies (C2) in Lemma \ref{lem:chenwang}.
  It then follows from Lemma \ref{lem:chenwang} that BSDE $(\xi,T,f_{n,q})$
  has a unique solution $(y^{n,q}_t,z^{n,q}_t)_{t\in\tT}$ in
  the space $\s^2\times \M^2$ for each fixed $n$, $q$.

  In the sequel, it follows from \ref{H:gMonotonicityInYWithU=0} on $f_n$
  and \eqref{eq:FnEnlargedByEAndU} 
  that $\pts$, for each $n,q\in\N$ and $y\in\rtn^k$,
  \begin{align*}
    \langle y,f_{n,q}(t,y)\rangle
    =\frac{q\vee|y|}{q}\langle \pi_q(y),f_n(t,\pi_q(y))-f_n(t,0)\rangle+
      \langle y,f_n(t,0)\rangle
    \leq |y||f_n(t,0)|\leq \big(K\me^{-t}+u(t)\phi(0)\big)|y|.
  \end{align*}
  Then assumption (A) is satisfied for the generator $f_{n,q}(t,y)$
  with $p=2$, $\mu(t)=\lambda(t)\equiv0$ and
  $f_t=K\me^{-t}+u(t)\phi(0)$.
  Thus, since $(y^{n,q}_t,z^{n,q}_t)_{t\in\tT}$ is the unique solution of BSDE
  $(\xi,T,f_{n,q})$ in $\s^2\times\M^2$, it follows from
  \eqref{eq:AssumptionFirstStep} and Proposition \ref{pro:estimate} with taking
  $r=t$ that there exists a universal positive constant $C_1$ such that for each $n,q\in\N$
  and $t\in\tT$,
  \begin{align}
    |y^{n,q}_t|^2+
    \EX
    \left[\left.
      \intT{t}|z^{n,q}_s|^2\dif s\right|\F_t
    \right]
    &\leq C_1\EX
      \left[\left.|\xi|^2+
        \left(
           \intT{t}[K\me^{-s}+\phi(0)u(s)]\dif s
        \right)^2\right|\F_t
      \right]\nonumber\\
    &\leq C_1\left[3K^2+2\phi^2(0)\left(\intT{0}u(s)\dif s\right)^2\right]:=a^2.
    \label{eq:YnqZnqBoundedByA}
  \end{align}
  Consequently, for any $q>a$, $(y^{n,q}_t,z^{n,q}_t)_{t\in\tT}$ does not
  depend on $q$. We then denote it by $(y^n_t,z^n_t)_{t\in\tT}$ and it is a solution
  of the following BSDE:
  \begin{equation}\label{eq:BSDE-fn-weak-convergence-before}
    y^n_t=\xi+\intT{t}f_n(s,y^n_s)\dif s-\intT{t}z^n_s\dif B_s,\quad t\in\tT.
  \end{equation}

  Furthermore, by \eqref{eq:YnqZnqBoundedByA} we have, for each $n\in\N$,
  \begin{equation}\label{eq:YnBoundedByAZnBoundedByA}
    \pts,\quad |y^n_t|^2\leq a^2 \quad\text{and}\quad
    \EX
    \left[
      \intT{0}|z^n_t|^2\dif t
    \right] \leq a^2.
  \end{equation}
  Assumption \ref{H:gGeneralGrowthInY} on $f_n$ and \eqref{eq:FnEnlargedByEAndU}
  yield that for each $n\in\N$,
  \[
    |f_n(t,y^n_t)|\leq |f_n(t,0)|+u(t)\phi(a)\leq K\me^{-t}+u(t)\big(\phi(0)+\phi(a)\big).
  \]
  Thus, we know that
  \begin{equation}\label{eq:SupYnFnZnInfite}
    \sup_n\EX
    \left[
      \sup_{t\in\tT}|y^n_t|^2+
      \left(
        \intT{0}|f_n(t,y^n_t)|\dif t
      \right)^2+
      \intT{0}|z^n_t|^2\dif t
    \right]<+\infty.
  \end{equation}
  Set $U^n_t:=f_n(t,y^n_t)$ for each $t\in\tT$. By \eqref{eq:SupYnFnZnInfite}
  we can conclude that there exists a subsequence of the sequence
  $\{(y^n_t,U^n_t,z^n_t)_{t\in\tT}\}^{\infty}_{n=1}$, still denoted by
  $\{(y^n_t,U^n_t,z^n_t)_{t\in\tT}\}^{\infty}_{n=1}$, that
  converges weakly in $\s^2(0,T;\rtn^k)\times\mathcal{L}^2(0,T;\rtn^k)\times\M^2(0,T;\rtn^{k\times d})$
  to a limit $(y_t,U_t,z_t)_{t\in\tT}$, where $\mathcal{L}^2(0,T;\rtn^k)$ denotes the set
  of $(\F_t)$-progressively measurable $\rtn^k$-valued processes
  $(U_t)_{t\in\tT}$ such that
  \begin{equation*}
    \|U\|_{\mathcal L^2}:=
    \left\{
      \EX\left[\left(\intT{0}|U_t|\dif t\right)^2\right]
    \right\}^{1\over 2}<+\infty.
  \end{equation*}
  In view of \eqref{eq:YnBoundedByAZnBoundedByA}, we have
  \begin{equation}\label{eq:YAndEXYInfinite}
    |y_t|\leq a, \quad \pts\quad \text{and} \quad
    \EX
    \left[
      \sup_{t\in\tT}|y_t|^2
    \right]<+\infty.
  \end{equation}
  In the sequel, take any bounded
  linear functional $\Phi(\cdot)$ defined on $L^2(\Omega,\F_T,\PR;\rtn^k)$.
  Then there exists a positive constant $b$ such that for each
  $(\overline{y}_s,\overline{U}_s,\overline{z}_s)_{s\in\tT}\in\s^2(0,T;\rtn^k)
  \times\mathcal L^2(0,T;\rtn^k)\times\M^2(0,T;\rtn^{k\times d})$ and $t\in\tT$,
  the following three inequalities hold true:
  \begin{equation*}
    |\Phi(\overline y _t)|\leq b\|\overline y _t\|_{L^2}\leq b\|\overline y\|_{\s^2},
  \end{equation*}
  \begin{equation*}
  \textstyle
    \left|
      \Phi
      \left(
        \intT{t}\overline{U}_s\dif s
      \right)
    \right|\leq b
    \left\|
      \intT{t}\overline U_s\dif s
    \right\|_{L^2}
      \leq b
      \left\|
        \overline{U}
      \right\|_{\mathcal L^2},
  \end{equation*}
  \begin{equation*}
  \textstyle
    \left|
      \Phi
      \left(
        \intT{t}\overline{z}_s\dif B_s
      \right)
    \right|\leq b
    \left\|
      \intT{t}\overline{z}_s\dif B_s
    \right\|_{L^2}
      \leq b\|\overline{z}\|_{\M^2}.
  \end{equation*}
  This means that for each $t\in\tT$,
  \[\Phi(\cdot_t),\quad \Phi(\textstyle\intT{t}\cdot\dif s),\quad \Phi(\textstyle\intT{t}\cdot\dif B_s)\]
  are bounded linear functionals defined respectively on
  $\s^2(0,T;\rtn^k)$, $\mathcal L^2(0,T;\rtn^k)$ and $\M^2(0,T;\rtn^{k\times d})$.
  Consequently, in view of the fact that $\{(y^n_t,U^n_t,z^n_t)_{t\in\tT}\}^\infty_{n=1}$
  converges weakly in $\s^2\times\mathcal L^2\times\M^2$ to the process $(y_t,U_t,z_t)_{t\in\tT}$,
  we have that for each $t\in\tT$,
  \begin{equation*}
    \lim_{n\to\infty}\Phi(y^n_t)=\Phi(y_t),
    \quad
    \lim_{n\to\infty}\Phi
    \left(
      \textstyle\intT{t}U^n_s\dif s
    \right)=\Phi
    \left(
      \textstyle\intT{t}U_s\dif s
    \right),
    \quad
    \lim_{n\to\infty}\Phi
    \left(
      \textstyle\intT{t}z^n_s\dif B_s
    \right)=\Phi
    \left(
      \textstyle\intT{t}z_s\dif B_s
    \right).
  \end{equation*}
  That is, for each $t\in\tT$, in the space $L^2(\Omega,\F_T,\PR;\rtn^k)$,
  $y^n_t$, $\intT{t}U_s^n\dif s$ and
  $\intT{t}z^n_s\dif B_s$ converge weakly to $y_t$, $\intT{t}U_s\dif s$ and
  $\intT{t}z_s\dif B_s$ respectively.
  Thus, taking weak limit in $L^2(\Omega,\F_T,\PR;\rtn^k)$ for BSDE
  \eqref{eq:BSDE-fn-weak-convergence-before} yields that for each $t\in\tT$,
  \begin{equation*}
    y_t=\xi+\intT{t}U_s\dif s-\intT{t}z_s\dif B_s, \quad\prs.
  \end{equation*}
  Then, noticing that $(y_t)_{t\in\tT}\in\s^2(0,T;\rtn^k)$ and
  the process $(\xi+\intT{t}U_s\dif s+\intT{t}z_s\dif B_s)_{t\in\tT}$
  is also continuous, we have, $\prs$,
  \begin{equation*}
    y_t=\xi+\intT{t}U_s\dif s-\intT{t}z_s\dif B_s, \quad t\in\tT.
  \end{equation*}

  Finally, by the following Lemma \ref{lem:AppendixWeakConvergenceU=f} we
  complete the proof of the first step.
  \begin{lem}\label{lem:AppendixWeakConvergenceU=f}
    $\pts$, $U_t=f(t,y_t)=g(t,y_t,V_t)$.
  \end{lem}
  The proof of Lemma \ref{lem:AppendixWeakConvergenceU=f}
  will be given in \ref{sec:Appendix}.

  {\medskip\noindent\bf Second step: }In this step we will prove that provided
  that $\xi\in L^p(\Omega,\F_T,\PR;\rtn^k)$,
  $V\in\M^p(0,T;\rtn^{k\times d})$ and that \ref{H:gContinuousInY},
  \ref{H:gGeneralizedGeneralGrowthInY}, \ref{H:gMonotonicityInYWithU=0},
  \ref{H:gLipschitzInZ} and \eqref{eq:AssumptionFirstStep} hold true,
  BSDE \eqref{eq:BSDE-g-yt-Vt} has a solution in $\s^2\times\M^2$.

  Assume now that $\xi\in L^p(\Omega,\F_T,\PR;\rtn^k)$,
  $V\in\M^p(0,T;\rtn^{k\times d})$ and that \ref{H:gContinuousInY},
  \ref{H:gGeneralizedGeneralGrowthInY}, \ref{H:gMonotonicityInYWithU=0},
  \ref{H:gLipschitzInZ} and \eqref{eq:AssumptionFirstStep} hold true.
  For some positive real $r'>0$, which will be chosen later, let $\theta_{r'}$ be a smooth
  function such that $0\leq \theta_{r'}(y)\leq 1$, $\theta_{r'}(y)=1$ for $|y|\leq r'$
  and $\theta_{r'}(y)=0$ as soon as $|y|\geq r'+1$. Now we define for each
  $(\omega,t,y)\in\Omega\times\tT\times\rtn^k$,
  \begin{equation*}
    h_n(t,y,V_t):=
    \theta_{r'}(y)\big(g(t,y,\pi_{n\me^{-t}}(V_t))-g(t,0,\pi_{n\me^{-t}}(V_t))\big)
    \frac{n\me^{-t}}{\psi_{r'+1}(t)\vee(n\me^{-t})}+g(t,0,V_t).
  \end{equation*}
  Then we have the following Lemma \ref{lem:AppendixHnProperty}, whose proof
  will be given in \ref{sec:Appendix}.
  \begin{lem}\label{lem:AppendixHnProperty}
    $h_n$ satisfies assumptions \ref{H:gContinuousInY}, \ref{H:gGeneralGrowthInY},
    \ref{H:gMonotonicityInY} and $|h_n(t,0,V_t)|\leq K\me^{-t}$.
  \end{lem}
  By Lemma \ref{lem:AppendixHnProperty} and Remark \ref{rmk:ChangeOfVariables},
  it follows from the first step that
  BSDE $(\xi,T,h_n)$ has a solution $(y^n_t,z^n_t)_{t\in\tT}$
  in the space $\s^2\times \M^2$. Thanks to assumption
  \ref{H:gMonotonicityInYWithU=0} on $g$ and \eqref{eq:AssumptionFirstStep},
  we get that for each $y\in\rtn^k$,
  \begin{equation*}
    \langle y,h_n(t,y,V_t)\rangle\leq \langle y,g(t,0,V_t)\rangle\leq|y||g(t,0,V_t)|\leq K\me^{-t}|y|,
  \end{equation*}
  which means that assumption (A) holds true for the generator $h_n(t,y,V_t)$
  with $\mu(t)=\lambda(t)\equiv 0$ and $f_t=K\me^{-t}$.
  Then it follows from Proposition \ref{pro:estimate} with $p=2$ that
  there exists a universal positive constant $C_2$ such that
  for each $n\in\N$ and $0\leq r\leq t\leq T$,
  \begin{equation*}
    \EX\left[\left.|y^n_t|^2\right|\F_r\right]
    +\EX\left[\left.
        \intT{t}|z^n_s|^2\dif s
      \right|\F_r
    \right]
   \leq C_2K^2:=r'^2.
  \end{equation*}
  Then we know that for each $n\in\N$,
  \begin{equation}\label{eq:YnAndZnBoundedByR}
    \pts,\quad |y^n_t|\leq r'\quad \text{and}\quad
    \EX
    \left[
      \intT{0}|z^n_s|^2\dif s
    \right]\leq r'^2.
  \end{equation}
  Hence, $(y^n_t,z^n_t)_{t\in\tT}$ is a solution of BSDE $(\xi,T,h'_n)$
  where
  \begin{equation*}
    h'_n(t,y,V_t):=\big(g(t,y,\pi_{n\me^{-t}}(V_t))-g(t,0,\pi_{n\me^{-t}}(V_t))\big)
    \frac{n\me^{-t}}{\psi_{r+1}(t)\vee(n\me^{-t})}+g(t,0,V_t).
  \end{equation*}
  It is clear that $h'_n$ satisfies \ref{H:gMonotonicityInYWithU=0} since
  $g$ satisfies it.

  In the sequel, for each $i,n\in\N$, we set $\hat{y}^{n,i}_t:=y^{n+i}_t-y^{n}_t$ and
  $\hat{z}^{n,i}_t:=z^{n+i}_t-z^n_t$. \Ito{} formula and assumption
  \ref{H:gMonotonicityInYWithU=0} on $h'_{n+i}$ yield that for each $t\in\tT$,
  \begin{equation*}
    |\hat{y}^{n,i}_t|^2+\intT{t}|\hat{z}^{n,i}_s|^2\dif s
    \leq 2\intT{t}|\hat{y}^{n,i}_s||h'_{n+i}(s,y^n_s,V_s)-h'_n(s,y^n_s,V_s)|\dif s
    -2\intT{t}\langle \hat{y}^{n,i}_s,\hat{z}^{n,i}_s\dif B_s\rangle,
  \end{equation*}
  from which it follows that
  \begin{equation}\label{eq:EnlargeForZAfterIto}
    \EX\left[\intT{0}|\hat{z}^{n,i}_s|^2\dif s\right]
    \leq 2\EX
         \left[
           \intT{0}|\hat{y}^{n,i}_s||h'_{n+i}(s,y^n_s,V_s)-h'_n(s,y^n_s,V_s)|\dif s
         \right],
  \end{equation}
  and
  \begin{align}\label{eq:EnlargeForYBeforeBDG}
    \EX\!
    \left[
      \sup_{t\in\tT}\!|\hat{y}^{n,i}_t|^2
    \right]\!\!
    \leq \!2\EX\!
    \left[
      \intT{0}\!\!|\hat{y}^{n,i}_s||h'_{n+i}(s,y^n_s,V_s)-h'_n(s,y^n_s,V_s)|\!\dif s
    \right]\!\!
    +2\EX\!
    \left[
      \sup_{t\in\tT}\!
      \left|
        \intT{t}\!\!\langle \hat y^{n,i}_s,\hat{z}^{n,i}_s\rangle\!\dif B_s
      \right|
    \right]\!\!.
  \end{align}
  Moreover, the Burkholder-Davis-Gundy (BDG for short in the remaining)
  inequality yields the existence of a constant $k$ such that
  \begin{align*}
    & 2\EX
    \left[
      \sup_{t\in\tT}
      \left|
        \intT{t}\langle \hat y^{n,i}_s,\hat{z}^{n,i}_s\rangle\dif B_s
      \right|
    \right]
    \leq 2k\EX
    \left[
      \left(
        \intT{0}|\hat y^{n,i}_s|^2|\hat{z}^{n,i}_s|^2\dif s
      \right)^{1\over 2}
    \right]\\
    & \leq 2k \EX
    \left[
      \sup_{s\in\tT}|\hat y^{n,i}_s|\cdot
      \left(
        \intT{0}|\hat z^{n,i}_s|^2\dif s
      \right)^{1\over 2}
    \right]
    \leq \frac{1}{2}\EX
    \left[
      \sup_{s\in\tT}|\hat y^{n,i}_s|^2
    \right]+2k^2\EX
    \left[
      \intT{0}|\hat{z}^{n,i}_s|^2\dif s
    \right].
  \end{align*}
  Putting the previous inequality into \eqref{eq:EnlargeForYBeforeBDG} we get
  \begin{equation}\label{eq:EnlargeForYAfterIto}
    \EX\left[\sup_{s\in\tT}|\hat{y}^{n,i}_s|^2\right]
    \leq 4\EX
         \left[
           \intT{0}|\hat{y}^{n,i}_s||h'_{n+i}(s,y^n_s,V_s)-h'_n(s,y^n_s,V_s)|\dif s
         \right]
         +4k^2\EX\left[\intT{0}|\hat{z}^{n,i}_s|^2\dif s\right].
  \end{equation}
  Combining with \eqref{eq:EnlargeForYAfterIto} and \eqref{eq:EnlargeForZAfterIto}
  and noticing by \eqref{eq:YnAndZnBoundedByR} that $\pts$,
  $|\hat{y}^{n,i}_t|\leq 2r'$, we get that there exists a constant $\overline k>0$
  such that
  \begin{equation}\label{eq:YAndZAndH}
    \EX\left[\sup_{s\in\tT}|\hat{y}^{n,i}_s|^2+\intT{0}|\hat{z}^{n,i}_s|^2\dif s\right]
    \leq r'\overline k\EX
         \left[
           \intT{0}|h'_{n+i}(s,y^n_s,V_s)-h'_n(s,y^n_s,V_s)|\dif s
         \right].
  \end{equation}
  Furthermore, we have the following Lemma \ref{lem:AppendixHnEnlargeForLebesgueConvrgence},
  whose proof will be provided in \ref{sec:Appendix}.
  \begin{lem}\label{lem:AppendixHnEnlargeForLebesgueConvrgence}
    For each $n,i\in\N$, the following inequality holds true:
    \begin{align*}
      |h'_{n+i}(s,y^n_s,V_s)-h'_n(s,y^n_s,V_s)|
      \leq 2v(s)|V_s|\one{|V_s|>n\me^{-s}}+2v(s)|V_s|\one{\psi_{r'+1}(s)>n\me^{-s}}
       +\psi_{r'+1}(s)\one{\psi_{r'+1}(s)>n\me^{-s}}.
    \end{align*}
  \end{lem}
  In view of $V\in\M^p(0,T;\rtn^{k\times d})$, $\psi_{r'+1}(t)\in L^1(\tT\times\Omega)$,
  $\intT{0}v^2(s)\dif s<+\infty$ and H\"older's inequality, Lebesgue's dominated
  convergence theorem yields that the right term of \eqref{eq:YAndZAndH} converges to
  $0$ as $n\to +\infty$. So $\{(y^n_t,z^n_t)_{t\in\tT}\}_{n=1}^\infty$
  is a Cauchy sequence in the space $\s^2\times \M^2$. Finally, passing
  to the limit in both sides of BSDE $(\xi,T,h'_n)$ under ucp (uniformly convergence
  in probability) yields a desired solution of BSDE \eqref{eq:BSDE-g-yt-Vt} in
  $\s^2\times\M^2$. The second step is then completed.

  {\medskip\noindent\bf Third step: }
  In this step we will eliminate the condition \eqref{eq:AssumptionFirstStep} used
  in the second step. 
  Under assumptions \ref{H:p-integrable} -- \ref{H:gGeneralizedGeneralGrowthInY},
  \ref{H:gMonotonicityInYWithU=0} and \ref{H:gLipschitzInZ}, we first define for each $n\in\N$,
  \begin{equation*}
    \xi^n:=\pi_n(\xi), \quad g^n(t,y,V_t):=g(t,y,V_t)-g(t,0,V_t)+\pi_{n\me^{-t}}(g(t,0,V_t)).
  \end{equation*}
  Then we can deduce that $|\xi^n|\leq n$, $|g^n(t,0,V_t)|\leq n\me^{-t}$
  and assumptions \ref{H:gContinuousInY}, \ref{H:gGeneralizedGeneralGrowthInY},
  \ref{H:gMonotonicityInYWithU=0} hold true for $g^n$. It follows from the
  second step that the following BSDE has a solution $(y^n_t,z^n_t)_{t\in\tT}$ in
  $\s^2\times \M^2$ such that $(y^n_t)_{t\in\tT}$ is a bounded process:
  \begin{equation}\label{eq:BSDEYntLp}
    y^n_t=\xi^n+\intT{t}g^n(s,y^n_s,V_s)\dif s-\intT{t}z^n_s\dif B_s, \quad t\in\tT.
  \end{equation}
  Obviously, $(y^n_t)_{t\in\tT}\in\s^p$. It follows from
  assumption \ref{H:gMonotonicityInYWithU=0} on $g$ that for each $y\in\rtn^k$,
  \begin{equation*}
    \langle y,g^n(t,y,V_t)\rangle\leq |y||\pi_{n\me^{-t}}(g(t,0,V_t))|\leq n\me^{-t}|y|,
  \end{equation*}
  which means that the generator $g^n$ satisfies assumption (A) with
  $\mu(t)=\lambda(t)\equiv 0$, $f_t=n\me^{-t}$. It then follows from
  Lemma \ref{lem:estimate-z}
  that $(z^n_t)_{t\in\tT}$ belongs to $\M^p$. Next we will prove the
  sequence $\{(y^n_t,z^n_t)_{t\in\tT}\}_{n=1}^{\infty}$ is a Cauchy sequence
  in the space $\s^p\times \M^p$.

  For each $n,m\in\N$,
  let $\hat\xi^{n,m}:=\xi^n-\xi^m$, $\hat{y}^{n,m}:=y^n-y^m$ and
  $\hat{z}^{n,m}:=z^n-z^m$. Then
  \begin{equation*}
    \hat{y}^{n,m}_t=\hat\xi^{n,m}+\intT{t}\hat{g}^{n,m}(s,\hat{y}^{n,m}_s,V_s)\dif s
    -\intT{t}\hat{z}^{n,m}_s\dif B_s,\quad t\in\tT,
  \end{equation*}
  where the generator
  $\hat{g}^{n,m}(t,y,V_t):=g^n(t,y+y^m_t,V_t)-g^m(t,y^m_t,V_t)$ for each $y\in \rtn^k$.
  In view of assumption \ref{H:gMonotonicityInYWithU=0} on $g^n$, we know that for each
  $y\in\rtn^k$ and $t\in\tT$,
  \begin{equation*}
    \langle y,\hat{g}^{n,m}(t,y,V_t)\rangle
    \leq \langle y,g^n(t,y^m_t,V_t)-g^m(t,y^m_t,V_t)\rangle
    \leq |y||\pi_{n\me^{-t}}(g(t,0,V_t))-\pi_{m\me^{-t}}(g(t,0,V_t))|.
  \end{equation*}
  Thus, assumption (A) is satisfied for the generator $\hat g^{n,m}$ with
  $u(t)=v(t)\equiv 0$ and
  $f_t=|\pi_{n\me^{-t}}(g(t,0,V_t))-\pi_{m\me^{-t}}(g(t,0,V_t))|$. Therefore,
  it follows from Proposition \ref{pro:estimate} with $r=t=0$ that there exists
  a positive constant $C^3_p$ depending only on $p$ such that for each $n,m\in\N$,
  \begin{align}
   & \EX
    \left[
      \sup_{s\in\tT}|\hat{y}^{n,m}_s|^p+
      \left(
        \intT{0}|\hat{z}^{n,m}_s|^2\dif s
      \right)^{p\over 2}
    \right]\nonumber\\
   &\leq C^3_p\EX\left[|\hat\xi^{n,m}|^p\right]
   +C^3_p\EX
    \left[
      \left(
        \intT{0}|\pi_{n\me^{-s}}(g(s,0,V_s))-\pi_{m\me^{-s}}(g(s,0,V_s))|\dif s
      \right)^p
    \right].\label{eq:YnmAndZnmLebesgueBefore}
  \end{align}
  Note that $\EX\left[|\xi^n-\xi|^p\right]\to 0$ as $n\to +\infty$ by Lebesgue's
  dominated convergence theorem. By assumption \ref{H:gLipschitzInZ}, we have that
  $|g(t,0,V_t)|\leq |g(t,0,0)|+v(t)|V_t|$, then H\"{o}lder's inequality yields
  that
  \begin{align*}
    \EX
    \left[
      \left(
        \intT{0}|g(t,0,V_t)|\dif t
      \right)^p
    \right]
    &\leq 2^{p-1}\EX
      \left[
        \left(
          \intT{0}|g(t,0,0)|\dif t
        \right)^p
      \right]\\
    &\quad\ +2^{p-1}
      \left(
        \intT{0}v^2(t)\dif t
      \right)^{p\over 2}\EX
      \left[
        \left(
          \intT{0}|V_t|^2\dif t
        \right)^{p\over 2}
      \right]<+\infty.
  \end{align*}
  Thus, by using Lebesgue's dominated convergence theorem once again, we get
  that as $n\to +\infty$,
  \begin{equation*}
    \EX\left[\left(\intT{0}|\pi_{n\me^{-t}}(g(t,0,V_t))-g(t,0,V_t)|\dif t\right)^p\right]\to 0.
  \end{equation*}
  Consequently, the right terms of \eqref{eq:YnmAndZnmLebesgueBefore} converges
  to $0$ as $n,m\to+\infty$.
  Hence, $\{(y^n_t,z^n_t)_{t\in\tT}\}_{n=1}^\infty$ is a Cauchy sequence
  in the space $\s^p\times \M^p$. Finally, by passing to the limit in both sides
  of \eqref{eq:BSDEYntLp} under ucp we can obtain a desired solution
  of BSDE \eqref{eq:BSDE-g-yt-Vt} in $\s^p\times\M^p$ .

  {\medskip\noindent\bf Fourth step: }
  In this step, we will finally complete the proof of the existence part of Theorem \ref{thm:MainResultLpSolution}.
  Assume now that $\xi\in L^p(\Omega,\F_T,\PR;\rtn^k)$ and that
  \ref{H:p-integrable} -- \ref{H:gGeneralizedGeneralGrowthInY}, \ref{H:gMonotonicityInYWithU=0}
  and \ref{H:gLipschitzInZ} hold true. It follows from the third step that
  for each $(V_t)_{t\in\tT}\in\M^p$,
  BSDE \eqref{eq:BSDE-g-yt-Vt} has a solution $(y_t,z_t)_{t\in\tT}$
  in $\s^p\times\M^p$. Take any $(U_t)_{t\in\tT}\in\s^p$
  and consider the following operator
  $
    \Psi:(U,V)\in \s^p\times\M^p\mapsto (y,z)
  $
  defined by
  \begin{equation*}
    y_t=\xi+\intT{t}g(s,y_s,V_s)\dif s-\intT{t}z_s\dif B_s,\quad t\in\tT.
  \end{equation*}
  Thus, we have constructed a mapping $\Psi$ from $\s^p\times \M^p$ to $\s^p\times \M^p$.
  Take another $(U',V')$ in the space $\s^p\times \M^p$ and set
  $(y',z'):=\Psi(U',V')$.
  Let us set $(\overline{U},\overline{V}):=(U-U',V-V')$ and
  $(\overline{y},\overline{z}):=(y-y',z-z')$ for notational convenience.
  Then $(\overline{y}_t,\overline{z}_t)_{t\in\tT}$
  is a solution of the following BSDE:
  \begin{equation}\label{eq:BSDE-contraction-estimate}
    \overline{y}_t=\intT{t}\overline{g}(s,\overline{y}_s)\dif s
    +\intT{t}\overline{z}_s\dif B_s,\quad t\in\tT,
  \end{equation}
  where the generator $\overline{g}(t,y):=g(t,y+y'_t,V_t)-g(t,y'_t,V'_t)$ for each
  $y\in\rtn^k$. It follows from assumptions \ref{H:gMonotonicityInYWithU=0} and
  \ref{H:gLipschitzInZ} on $g$ that
  \begin{align*}
    \langle y,\overline{g}(t,y)\rangle
     \leq \langle y, g(t,y'_t,V_t)-g(t,y'_t,V'_t)\rangle
     \leq v(t)|\overline{V}_t||y|,
  \end{align*}
  which means that $\overline{g}$ satisfies assumption (A) with
  $\mu(t)=\lambda(t)\equiv 0$ and $f_t=v(t)|\overline{V}_t|$. Thus,
  it follows from Proposition \ref{pro:estimate} and H\"{o}lder's inequality that
  there exists a positive constant $C^4_p$ depending only on $p$ such that
  \begin{align*}
    \EX\!\!
    \left[
      \sup_{t\in\tT}\!|\overline{y}_t|^p+\!
      \left(
        \intT{0}\!\!|\overline{z}_t|^2\dif t
      \right)^{p\over 2}
    \right]\!\!
    \leq C^4_p\EX\!\!
     \left[
       \left(
         \intT{0}v(t)|\overline{V}_t|\dif t
       \right)^p
     \right]\!\!
    \leq C^4_p\left(\intT{0}v^2(t)\dif t\right)^{p\over 2}\!
          \EX\!\!\left[\left(\intT{0}|\overline{V}_t|^2\dif t\right)^{p\over 2}\right]\!.
  \end{align*}
  Hence, if $\delta:=C^4_p\big(\intT{0}v^2(t)\dif t\big)^{p/2}<1$, we have the
  strict contraction in the space $\s^p\times \M^p$ as follows:
  \begin{equation*}
    \|(\overline{y},\overline{z})\|^p_{\s^p\times\M^p}
    <\delta\|(\overline{U},\overline{V})\|^p_{\s^p\times\M^p}.
  \end{equation*}
  Then by the fixed point theorem BSDE \eqref{eq:BSDEs}
  has a unique solution in $\s^p\times\M^p$.

  In the general case about $\delta$, in view of
  the fact that $\intT{0}v^2(t)\dif t<+\infty$ and Proposition \ref{pro:estimate},
  we can follow exactly the proof procedure of the existence part in
  \citet{FanJiang2011Stochastics}. That is, we can subdivide the
  interval $\tT$ into a finite number of small intervals like
  $[T_i,T_{i+1}]$ $(i=0,1,\cdots,N-1)$ such that $0=T_0<T_1<\cdots<T_{N-1}<T_N=+\infty$ and
  for each $i=0,1,\cdots,N-1$,
  \[C^4_p\left(\intT[T_{i+1}]{T_i}v^2(t)\dif t\right)^{p\over 2}<1.\]
  And in every small interval there exists a strict contraction in
  the space $\s^p\times \M^p$. Then we complete the proof of the existence part
  of Theorem \ref{thm:MainResultLpSolution}.
\end{proof}

\begin{rmk}
  Motivated by Remark 9.3 in \citet{Touzi2013StochasticControlNote} we can also
  consider Picard's iteration procedure to show the fourth step.
  Set $(y^0,z^0):=(0,0)$ and define $\{(y^n_t,z^n_t)_{t\in\tT}\}^\infty_{n=1}$
  recursively, in view of the third step, for each $n\geq 0$,
  \begin{equation*}
    y^{n+1}_t=\xi+\intT{t}g(s,y^{n+1},z^n_s)\dif s-\intT{t}z^{n+1}_s\dif B_s,\quad t\in\tT.
  \end{equation*}
  Set $\delta y^n:=y^{n+1}-y^n$, $\delta z^n:=z^{n+1}-z^n$, then we have
  \begin{equation*}
    \delta y^n_t=\intT{t}\delta g^n(s,\delta y^n_s)\dif s-\intT{t}\delta z^n_s\dif B_s,\quad t\in\tT,
  \end{equation*}
  where $\delta g^n(s,y):=g(s,y+y^n_s,z^n_s)-g(s,y^n_s,z^{n-1}_s)$ for each $y\in\rtn^k$.
  Since $\delta g^n$ satisfies assumption (A) with $\mu(t)=u(t)$, $\lambda(t)\equiv 0$
  and $f_t=v(t)|\delta z^{n-1}_t|$, it follows from Proposition \ref{pro:estimate}
  and an induction argument that
  \begin{equation*}
    \|(\delta y^n,\delta z^n)\|^p_{\s^p\times\M^p}
    \leq \left[C'_p\left(\intT{0}v^2(s)\dif s\right)^{p\over 2}\right]^{n-1}\|(\delta y^1,\delta z^1)\|^p_{\s^p\times\M^p}.
  \end{equation*}
  Then by subdividing the time interval $[0,T]$ into a finite number of small intervals
  like $[T_i,T_{i+1}]$ $(i=0,1,\cdots,N-1)$ such that $0=T_0<T_1<\cdots<T_{N-1}<T_N=+\infty$
  and for each $i=0,1,\cdots,N-1$, $C'_p\big(\intT[T_{i+1}]{T_i}v^2(s)\dif s\big)^{p/2}<1$,
  we can deduce that $\{(y^n_t,z^n_t)_{t\in[T_i,T_{i+1}]}\}_{n=1}^\infty$ is a
  Cauchy sequence in $\s^p(T_{i-1},T_i;\rtn^k)\times\M^p(T_{i-1},T_i;\rtn^{k\times d})$,
  and then the limit process $(y_t,z_t)_{t\in\tT}$ is a solution of BSDE \eqref{eq:BSDEs}
  in $\s^p\times\M^p$.
\end{rmk}

\begin{ex}\label{ex:LpUnboundedParameters}
  Let $0\leq T<+\infty$, $k=1$ and
  $$g(t,y,z)=|\ln t\,|(-\me^{y}+|y|)+\frac{|z|}{\sqrt[4]{t\,}}+|B_t|.$$
  It is not difficult to check that $g$ satisfies assumptions
  \ref{H:p-integrable} -- \ref{H:gLipschitzInZ} with $u(t)=|\ln t\,|$ and $v(t)=1/\sqrt[4]{t\,}$.
  Then by Theorem
  \ref{thm:MainResultLpSolution} we know that for each
  $\xi\in L^p(\Omega,\F_T,\PR;\rtn^k)$, BSDE $(\xi,T,g)$ has a unique
  solution in $\s^p\times\M^p$. But because $u(t)$ and $v(t)$ are unbounded,
  this conclusion can not be obtained by Theorem 4.2 in \citet{BriandDelyonHu2003SPA}.
\end{ex}

\begin{ex}
  Let $0\leq T\leq +\infty$, $k=2$ and
  \begin{equation*}
    g(t,y,z)=t^2\me^{-t}
    \begin{bmatrix}
      \displaystyle-y^3_1+y_2\\
      \displaystyle-y^5_2-y_1
    \end{bmatrix}
    +\frac{1}{\sqrt{1+t^2}}
    \begin{bmatrix}
      \displaystyle|z_1|\\
      \displaystyle|z_2|
    \end{bmatrix}
    +\frac{t^2}{t^4+1}
    \begin{bmatrix}
      \displaystyle 1\\
      \displaystyle 1
    \end{bmatrix},
  \end{equation*}
  where $y_i$ and $z_i$ $(i=1,2)$ stand for the $i$th component of the vector $y$
  and the $i$th row of the matrix $z$ respectively.
  It is not hard to verify that $g$ satisfies
  assumptions \ref{H:p-integrable} -- \ref{H:gLipschitzInZ}
  with $u(t)=t^2\me^{-t}$, $v(t)=1/\sqrt{1+t^2}$. Thus, by Theorem
  \ref{thm:MainResultLpSolution} we know that for each
  $\xi\in L^p(\Omega,\F_T,\PR;\rtn^k)$, BSDE $(\xi,T,g)$ has a unique
  solution in $\s^p\times\M^p$. It should be noted that this conclusion can
  not be obtained by Theorem 4.2 in
  \citet{BriandDelyonHu2003SPA} when $T=+\infty$.
\end{ex}

\section{$L^1$ solution}
\label{sec:IntegrableParameters}
In this section we will give an existence and uniqueness result for
$L^1$ solutions of BSDE \eqref{eq:BSDEs}. Here, we suppose that the generator $g$ is
monotonic and has a general growth in $y$, and is Lipschitz continuous
and has a kind of sublinear growth in $z$, which are both non-uniform with
respect to $t$.

We first introduce the following assumptions on the generator $g$,
where $0\leq T\leq +\infty$.
\begin{enumerate}
\renewcommand{\theenumi}{(H\arabic{enumi}')}
\renewcommand{\labelenumi}{\theenumi}

  \item\label{H:OnlyIntegrable}
          $\EX
          \left[
            \intT{0}|g(t,0,0)|\dif t
          \right]<+\infty$;
          \vspace{-5pt}
\renewcommand{\theenumi}{(H\arabic{enumi})}
\setcounter{enumi}{5}
  \item\label{H:gSublinearGrowthInZ}
        There exist an $\alpha\in(0,1)$ and a deterministic function
        $\gamma(t):\tT\mapsto \rtn^+$ with
        $\intT{0}\big(\gamma(t)+\gamma^{1/(1-\alpha)}(t)+\gamma^{2/(2-\alpha)}(t)\big)\dif t<+\infty$
        such that $\pts$, for each $y\in \rtn^k$ and $z\in\rtn^{k\times d}$,
        $$|g(t,y,z)-g(t,y,0)|\leq \gamma(t)(g_t+|y|+|z|)^\alpha,$$
        where $(g_t)_{t\in[0,T]}$ is a nonnegative and $(\F_t)$-progressively
        measurable process with $\EX[\int^T_0 g_t\dif t]<+\infty$.
\end{enumerate}

The following Theorem \ref{thm:MainResultL1Solution} is the main result of this
section.
\begin{thm}\label{thm:MainResultL1Solution}
  Let $0\leq T\leq +\infty$ and assumptions \ref{H:OnlyIntegrable},
  \ref{H:gContinuousInY} -- \ref{H:gSublinearGrowthInZ} on the generator $g$ hold.
  Then for each $\xi\in L^1(\Omega,\F_T,\PR;\rtn^k)$,
  BSDE \eqref{eq:BSDEs} has a solution $(y_t,z_t)_{t\in\tT}$
  in $\s^\beta\times\M^\beta$
  for $\beta\in(0,1)$ with $(y_t)_{t\in\tT}$ belonging to the class (D),
  which is unique in $\s^\beta\times\M^\beta$ for $\beta\in(\alpha,1)$.
\end{thm}

The proof of Theorem \ref{thm:MainResultL1Solution} is completed with the help
of Theorem \ref{thm:MainResultLpSolution} in the sprit of Theorems 6.2 and 6.3
in \citet{BriandDelyonHu2003SPA}.
In view of Remark \ref{rmk:ChangeOfVariables}, we shall always assume that $g$
satisfies \ref{H:OnlyIntegrable}, \ref{H:gContinuousInY}, \ref{H:gGeneralizedGeneralGrowthInY},
\ref{H:gMonotonicityInYWithU=0}, \ref{H:gLipschitzInZ} and \ref{H:gSublinearGrowthInZ}.
As usual, let us first show the uniqueness part.

\begin{proof}[\bf Proof of the uniqueness part of Theorem \ref{thm:MainResultL1Solution}]
  Assume that both $(y_t,z_t)_{t\in\tT}$ and $(y'_t,z'_t)_{t\in\tT}$ are
  solutions of BSDE \eqref{eq:BSDEs} such that both $(y_t)_{t\in\tT}$ and
  $(y'_t)_{t\in\tT}$ belong to the class (D), and both $(y_t,z_t)_{t\in\tT}$ and $(y'_t,z'_t)_{t\in\tT}$
  belong to $\s^\beta\times\M^\beta$ for some $\beta\in(\alpha, 1)$.
  Let us fix $n\in \N$ and denote $\tau_n$ the stopping time
  $$\tau_n=\inf\left\{t\in\tT:\intT[t]{0}(|z_s|^2+|z'_s|^2)\dif s\geq n\right\}\wedge T.$$

   Lemma \ref{lem:Tanaka-Briand} leads to the following
   inequality with setting $\hat{y}:=y-y'$ and
   $\hat{z}:=z-z'$,
   \begin{align}\label{eq:tanaka-multi-no-exp}
     |\hat{y}_{t\wedge\tau_n}|
     \leq
     &\ |\hat{y}_{\tau_n}|+\intT[\tau_n]{t\wedge\tau_n}
          |\hat{y}_s|^{-1}\one{|\hat{y}_s|\neq 0}
          \langle \hat{y}_s,g(s,y_s,z_s)-g(s,y'_s,z'_s)\rangle\dif s\nonumber\\
     &\ -\intT[\tau_n]{t\wedge\tau_n}|\hat{y}_s|^{-1}
        \one{|\hat{y}_s|\neq 0}\langle \hat{y}_s,\hat{z}_s\dif B_s\rangle.
   \end{align}
   We first enlarge the inner product including $g$ via assumptions
   \ref{H:gMonotonicityInYWithU=0} and \ref{H:gSublinearGrowthInZ} as follows:
   \begin{align*}
     |\hat{y}_s|^{-1}\one{|\hat{y}_s|\neq 0}
      \langle \hat{y}_s,g(s,y_s,z_s)-g(s,y'_s,z'_s)\rangle
     \leq |g(s,y'_s,z_s)-g(s,y'_s,z'_s)|
     \leq2\gamma(s)
      \left(
        g_s+|y'_s|+|z'_s|+|z_s|
      \right)^\alpha.
   \end{align*}
   Putting the previous inequality into \eqref{eq:tanaka-multi-no-exp}
   and then taking conditional expectation with respect to $\F_t$ in both sides,
   we have that for each $n\in\N$,
   \begin{equation*}
       |\hat{y}_{t\wedge\tau_n}|
      \leq \EX
       \left[
         |\hat{y}_{\tau_n}|+\left.2\intT{0}\gamma(s)
         \left(
           g_s+|y'_s|+|z'_s|+|z_s|
         \right)^\alpha\dif s\right|\F_t
       \right].
   \end{equation*}
   Now sending $n\to\infty$ and noticing that $\tau_n\to T$, $(\hat{y}_t)_{t\in\tT}$ belongs to the
   class (D) and $\prs$, $\hat{y}_T=0$, we know that for each $t\in\tT$,
   \begin{equation*}
     |\hat{y}_t|\leq 2\EX
     \left[
       \left.\intT{0}\gamma(s)\left(g_s+|y'_s|+|z'_s|+|z_s|\right)^\alpha\dif s\right|\F_t
     \right].
   \end{equation*}
   Moreover, noticing that $\beta/\alpha>1$, Doob's inequality yields
   that there exists a positive constant $C^\alpha_\beta$ depending on $\alpha$
   and $\beta$ such that
   \begin{equation}\label{eq:DoobResultHatYBetaAlpha}
     \EX
     \left[
       \sup_{t\in\tT}|\hat{y}_t|^{\beta\over\alpha}
     \right]\leq C^\alpha_\beta\EX
     \left[
       \left(
         \intT{0}\gamma(s)(g_s+|y'_s|+|z'_s|+|z_s|)^\alpha\dif s
       \right)^{\beta\over\alpha}
     \right].
   \end{equation}
   We set the random variable $G:=\intT{0}\gamma(s)(g_s+|y'_s|+|z'_s|+|z_s|)^\alpha\dif s$.
   Then $G\in L^{\beta/\alpha}(\Omega,\F_T,\PR;\rtn)$. Indeed, H\"{o}lder's
   inequality yields that
   \begin{equation*}
     \intT{0}\gamma(s)g^\alpha_s\dif s
     \leq
     \left(
       \intT{0}\gamma^{1\over 1-\alpha}(s)\dif s
     \right)^{1-\alpha}
       \left(
         \intT{0}g_s\dif s
       \right)^\alpha,
   \end{equation*}
   \begin{equation*}
     \intT{0}\gamma(s)|z_s|^\alpha\dif s
     \leq \left(\intT{0}\gamma^{2\over 2-\alpha}(s)\dif s\right)^{2-\alpha\over 2}
     \left(\intT{0}|z_s|^2\dif s\right)^{\alpha\over 2},
   \end{equation*}
   and $z'_s$ has the similar estimate. Besides,
   \begin{align*}
     \intT{0}\gamma(s)|y'_s|^\alpha\dif s
     \leq
     \intT{0}\gamma(s)\dif s
     \cdot\sup_{t\in\tT}|y'_t|^{\alpha}.
   \end{align*}
   Consequently, noticing that $\EX[\intT{0}g_s\dif s]<+\infty$,
   $(y'_t)_{t\in\tT}$ belongs to the space
   $\s^\beta$, $(z_t)_{t\in\tT}$ and $(z'_t)_{t\in\tT}$ belongs to the space
   $\M^\beta$, we have that $G\in L^{\beta/\alpha}$. It then follows from
   \eqref{eq:DoobResultHatYBetaAlpha} that $(\hat{y}_t)_{t\in\tT}$ belongs to
   the space $\s^{\beta/\alpha}$.
   Furthermore, note that $(\hat y_t, \hat z_t)_{t\in\tT}$ is a solution of
   the following BSDE:
   \begin{equation*}
     \hat y_t=\intT{t}\hat g(s,\hat y_s,\hat z_s)\dif s
     -\intT{t}\hat z_s\dif B_s,\quad t\in\tT,
   \end{equation*}
   where the generator $\hat g(t,\hat y,\hat z):=g(t,y+y'_t,z+z'_t)-g(t,y'_t,z'_t)$
   for each $(y,z)\in\rtn^k\times\rtn^{k\times d}$. It follows from
   assumptions \ref{H:gMonotonicityInYWithU=0} and \ref{H:gLipschitzInZ} on the
   generator $g$ that
   $$\langle y,\hat g(t,\hat y,\hat z)\rangle
   \leq \langle y,g(t,y'_t,z+z'_t)-g(t,y'_t,z'_t)\rangle\leq v(t)|y||z|,$$
   which means that assumption (A) is satisfied for the generator $\hat g$
   with $\mu(t)\equiv 0$,
   $\lambda(t)=v(t)$ and $f_t\equiv 0$. It then follows from
   Proposition \ref{pro:estimate} with $p=\beta/\alpha$ that
   $(\hat y_t,\hat z_t)_{t\in\tT}=(0,0)\in\s^{\beta/\alpha}\times \M^{\beta/\alpha}$.
   By now the uniqueness part is proved completely.
\end{proof}

\begin{proof}[\bf Proof of the existence part of Theorem \ref{thm:MainResultL1Solution}]
  The proof will be split into two steps. The first step deals with the case that
  the generator $g$ is independent of the variable $z$ under assumptions
  \ref{H:OnlyIntegrable}, \ref{H:gContinuousInY}, \ref{H:gGeneralizedGeneralGrowthInY}
  and \ref{H:gMonotonicityInYWithU=0}, and the
  second step considers the general case.

  {\medskip\noindent\bf First step:} We assume that $g$ is independent of
  $z$. For each $n\in \N$, we set
  \begin{equation*}
    \xi^n:=\pi_n(\xi),\quad
    g^n(t,y):=g(t,y)-g(t,0)+\pi_{n\me^{-t}}(g(t,0)).
  \end{equation*}
  Note that $|\xi^n|\leq n$ and $g^n(t,y)$ satisfies assumptions
  \ref{H:p-integrable} -- \ref{H:gGeneralizedGeneralGrowthInY} and
  \ref{H:gMonotonicityInYWithU=0}. It follows from
  Theorem \ref{thm:MainResultLpSolution} that the following BSDE
  \eqref{eq:BSDE-gn-independent-z} has a unique
  solution $(y^n_t,z^n_t)_{t\in\tT}$ in $\s^2\times \M^2$,
  \begin{equation}\label{eq:BSDE-gn-independent-z}
    y^n_t=\xi^n+\intT{t}g^n(s,y^n_s)\dif s-\intT{t}z_s^n\dif B_s,\quad t\in\tT.
  \end{equation}
  For each $i$, $n\in\N$, we set $\hat{y}^{n,i}:=y^{n+i}-y^n$,
  $\hat z^{n,i}:=z^{n+i}-z^n$ and $\hat\xi^{n,i}:=\xi^{n+i}-\xi^n$. It then follows
  from Lemma \ref{lem:Tanaka-Briand} with taking $p=1$ that for each $t\in\tT$,
  \begin{align}
    |\hat{y}^{n,i}_t|
    \leq
    &\ |\hat{\xi}^{n,i}|+\intT{t}|\hat y^{n,i}_s|^{-1}\one{|\hat y^{n,i}_s|\neq 0}
     \langle \hat y^{n,i}_s,g^{n+i}(s,y^{n+i}_s)-g^n(s,y^n_s)\rangle\dif s\nonumber\\
    &\ -\intT{t}|\hat y^{n,i}_s|^{-1}\one{|\hat y^{n,i}_s|\neq 0}\langle \hat y^{n,i}_s,\hat z^{n,i}_s\dif B_s\rangle.\label{eq:TanakaNoStoppingTime}
  \end{align}
  As before, the inner product including $g$ can be enlarged via assumption
  \ref{H:gMonotonicityInYWithU=0} on $g^{n+i}$ as follows:
  \begin{equation*}
    |\hat y^{n,i}_s|^{-1}\one{|\hat y^{n,i}_s|\neq 0}
    \langle \hat y^{n,i}_s,g^{n+i}(s,y^{n+i}_s)-g^n(s,y^n_s)\rangle
    \leq |g^{n+i}(s,y^n_s)-g^n(s,y^n_s)|.
  \end{equation*}
  Putting the previous inequality into \eqref{eq:TanakaNoStoppingTime}
  and taking conditional expectation with respect to $\F_t$ in both sides,
  we can get that for each $t\in\tT$,
  \begin{align*}
    |\hat y_t^{n,i}|
    &\leq \EX
    \left[
      |\hat{\xi}^{n,i}|+\left.\intT{0}|g^{n+i}(s,y^n_s)-g^n(s,y^n_s)|\dif s\right|\F_t
    \right]\\
    &\leq \EX
     \left[
       |\xi|\one{|\xi|>n}+\left.\intT{0}|g(s,0)|\one{|g(s,0)|>n\me^{-s}}\dif s\right|\F_t
     \right].
  \end{align*}
  Thus, we have
  \begin{equation*}
    \|\hat{y}^{n,i}_t\|_1\leq \EX
     \left[
       |\xi|\one{|\xi|>n}+\intT{0}|g(s,0)|\one{|g(s,0)|>n\me^{-s}}\dif s
     \right].
  \end{equation*}
  And according to Lemma 6.1 in \citet{BriandDelyonHu2003SPA} we know that for any
  $\beta\in(0,1)$,
  \begin{align*}
    \EX
    \left[
      \sup_{t\in\tT}|\hat{y}_t^{n,i}|^\beta
    \right]\leq\frac{1}{1-\beta}
    \left(\EX
      \left[
        |\xi|\one{|\xi|>n}+\intT{0}|g(s,0)|\one{|g(s,0)|>n\me^{-s}}\dif s
      \right]
    \right)^{\beta}.
  \end{align*}
  Therefore, note that $\EX[|\xi|+\intT{0}|g(s,0)|\dif s]<+\infty$, the process
  sequence $\{(y^n_t)_{t\in\tT}\}_{n=1}^{\infty}$ is a Cauchy sequence both
  under $\|\cdot\|_1$ and in the space $\s^\beta$. Let
  $(y_t)_{t\in\tT}$ denote the limit of this process sequence, then
  $(y_t)_{t\in\tT}$ belongs to the class (D) and
  $\s^\beta$ for any $\beta\in(0,1)$.

  Furthermore, note that $(\hat y_t^{n,i},\hat z_t^{n,i})_{t\in\tT}$ is a
  solution of the following BSDE:
  \begin{equation*}
    \hat{y}_t^{n,i}=\hat \xi^{n,i}+\intT{t}\hat{g}(s,\hat y_s^{n,i})\dif s-\intT{t}\hat{z}_s^{n,i}\dif B_s, \quad t\in\tT,
  \end{equation*}
  where the generator $\hat g(t,y):=g^{n+i}(t,y+y^n_t)-g^n(t,y^n_t)$ for each
  $y\in\rtn^k$. It follows from assumption \ref{H:gMonotonicityInYWithU=0} on $g^{n+i}$
  that
  \begin{align*}
    \langle y,\hat g(t,y)\rangle
    \leq \langle y,g^{n+i}(t,y^n_t)-g^n(t,y^n_t)\rangle
    \leq |y||g(t,0)|\one{|g(t,0)|>n\me^{-t}},
  \end{align*}
  which means that $\hat g$ satisfies assumption (A) with
  $\mu(t)=\lambda(t)\equiv 0$ and $f_t=|g(t,0)|\one{|g(t,0)|>n\me^{-t}}$.
  It then follows from Lemma \ref{lem:estimate-z} with $p=\beta$ that for any $\beta\in(0,1)$,
  there exists a constant $C^1_\beta$ depending on $\beta$ such that
  \begin{equation*}
    \EX
    \left[
      \left(
        \intT{0}|\hat z_t^{n,i}|^2\dif t
      \right)^{\beta\over 2}
    \right]\leq C^1_\beta\EX
    \left[
      \sup_{t\in\tT}|\hat{y}_t^{n,i}|^\beta+
      \left(
        \intT{0}|g(t,0)|\one{|g(t,0)|>n\me^{-t}}\dif t
      \right)^\beta
    \right].
  \end{equation*}
  In view of the previous inequality, we know that for any $\beta\in(0,1)$,
  the process sequence
  $\{(z^n_t)_{t\in\tT}\}_{n=1}^{\infty}$ is a Cauchy sequence in the space
  $\M^\beta$. Let $(z_t)_{t\in\tT}$ denote the limit which belongs to the
  space $\M^\beta$ for any $\beta\in(0,1)$, then $\intT{t}z^n_s\dif B_s$
  converges to $\intT{t}z_s\dif B_s$ under ucp.

  Finally, since $y\mapsto g(t,y)$ is continuous, we can obtain that
  $(y_t,z_t)_{t\in\tT}$ is a desired solution
  of BSDE \eqref{eq:BSDEs} via taking the limit in
  both sides of BSDE \eqref{eq:BSDE-gn-independent-z} under ucp.

  {\medskip\noindent\bf Second step: } The general case, i.e., $g$ may depend on $z$.

  The next proof procedure will use Picard's iterative procedure.
  Let us set $(y^0,z^0):=(0,0)$. It is not hard to verify that for each
  $(z_t)_{t\in\tT}\in\M^\beta$ for any $\beta\in(0,1)$, the
  generator $g(t,y,z_t)$ satisfies \ref{H:OnlyIntegrable}, \ref{H:gContinuousInY},
  \ref{H:gGeneralizedGeneralGrowthInY} and \ref{H:gMonotonicityInYWithU=0}.
  Thus, we can define the process sequence
  $\{(y^n_t,z^n_t)_{t\in\tT}\}_{n=1}^\infty$ recursively, in view
  of the first step, for each $n\geq 0$,
  \begin{equation}\label{eq:BSDEYnPicardIterative}
    y^{n+1}_t=\xi+\intT{t}g(s,y^{n+1}_s,z^n_s)\dif s-\intT{t}z^{n+1}_s\dif B_s, \quad t\in\tT,
  \end{equation}
  where for each $n\geq 0$, $(y^n_t)_{t\in\tT}$ belongs to the class (D) and $(y^n_t,z^n_t)_{t\in\tT}$ belongs
  to $\s^\beta\times \M^\beta$ for any $\beta\in(0,1)$.

  For each $n\in\N$, arguing as in the proof of the uniqueness part of Theorem
  \ref{thm:MainResultL1Solution}, we can deduce, in view
  of assumptions \ref{H:gMonotonicityInYWithU=0} and \ref{H:gSublinearGrowthInZ},
  that for each $t\in\tT$,
  \begin{align}
    |y^{n+1}_t-y^n_t|
   &\leq
      \EX\left[\left.
           \intT{0}|g(s,y^n_s,z^n_s)-g(s,y^n_s,z^{n-1}_s)|\dif s\right|\F_t
         \right]\nonumber\\
   &\leq 2\EX
    \left[\left.
       \intT{0}\gamma(s)\left(g_s+|y^n_s|+|z^n_s|+|z^{n-1}_s|\right)^\alpha\dif s\right|\F_t
    \right].
    \label{eq:YnBoundedByAlphaWithConditionalEX}
  \end{align}
  We set the random variable
  \begin{equation*}
    I_n:=\intT{0}\gamma(s)\left(g_s+|y^n_s|+|z^n_s|+|z^{n-1}_s|\right)^\alpha\dif s.
  \end{equation*}
  Similar to the proof procedure of the uniqueness part of Theorem
  \ref{thm:MainResultL1Solution}, we can prove that $I_n$ belongs to
  $L^q(\Omega,\F_T,\PR;\rtn)$ as soon as $\alpha q<1$ with
  $q>1$. Furthermore, in view of Doob's inequality, for some $q>1$ such that
  $\alpha q<1$, by \eqref{eq:YnBoundedByAlphaWithConditionalEX}
  we can deduce that there exists a positive constant
  $c_q$ depending only on $q$ such that
  \begin{equation*}
    \EX
    \left[
      \sup_{t\in\tT}|y^{n+1}_t-y^n_t|^q
    \right]\leq c_q\EX\left[I^q_n\right]<+\infty,
  \end{equation*}
  which implies that $(y^{n+1}_t-y^n_t)_{t\in\tT}$ belongs to the space
  $\s^q$ for some $q>1$.

  In the sequel, for each $n\in\N$, we set $\hat{y}^n:=y^{n+1}-y^n$ and
  $\hat z^n:=z^{n+1}-z^n$,
  then $(\hat y^n_t,\hat z^n_t)_{t\in\tT}$ is a solution of the following BSDE:
  \begin{equation*}
    \hat y^n_t=\intT{t}g^n(s,\hat y^n_s)\dif s-\intT{t}\hat z^n_s\dif B_s, \quad t\in\tT,
  \end{equation*}
  where the generator $g^n(t,y):=g(t,y+y^n_t,z^n_t)-g(t,y^n_t,z^{n-1}_t)$ for
  each $y\in\rtn^k$.
  It follows from assumptions \ref{H:gMonotonicityInYWithU=0} and
  \ref{H:gSublinearGrowthInZ} that
  \begin{align*}
    \langle y,g^n(t,y)\rangle
    \leq |y||g(t,y^n_t,z^n_t)-g(t,y^n_t,z^{n-1}_t)|
    \leq 2|y|\gamma(t)\left(g_t+|y^n_t|+|z^n_t|+|z^{n-1}_t|\right)^\alpha,
  \end{align*}
  which means that the generator $g^n$ satisfies assumption (A) with
  $\mu(t)=\lambda(t)\equiv 0$,
  $f_t=2\gamma(t)(g_t+|y^n_t|+|z^n_t|+|z^{n-1}_t|)^\alpha$ and $p=q$
  since $I_n$ belongs to $L^q$.
  Thus, Lemma \ref{lem:estimate-z} yields that $(\hat z_t^n)_{t\in\tT}$ belongs to the
  space $\M^q$. Besides, in view of assumptions \ref{H:gMonotonicityInYWithU=0}
  and \ref{H:gLipschitzInZ} we have $\langle y,g^n(t,y)\rangle \leq v(t)|y||\hat z_t^{n-1}|$.
  Thus, Proposition \ref{pro:estimate} with $p=q$, $\mu(t)=\lambda(t)\equiv 0$
  and $f_t=v(t)|\hat z^{n-1}_t|$ yields that there exists
  a constant $C_q>0$ depending only on $q$ such that for each $n\geq 2$,
  \begin{align*}
    \|(\hat y^n,\hat z^n)\|^q_{\s^q\times\M^q}
    \leq C_q\EX
     \left[
       \left(
         \intT{0}v(t)|\hat z^{n-1}_t|\dif t
       \right)^q
     \right].
  \end{align*}
  Then by H\"{o}lder's inequality we get
  \begin{align*}
    \|(\hat y^n,\hat z^n)\|^q_{\s^q\times\M^q}
    &\leq C_q\left(\intT{0}v^2(t)\dif t\right)^{q\over 2}\EX
     \left[
       \left(
         \intT{0}|\hat z^{n-1}_t|^2\dif t
       \right)^{q\over 2}
     \right],
  \end{align*}
  from which it follows that for each $n\geq 2$,
  \begin{equation*}
    \|(\hat y^n,\hat z^n)\|^q_{\s^q\times\M^q}\leq
    \left[
      C_q\left(\intT{0}v^2(t)\dif t\right)^{q\over 2}
    \right]^{n-1}\|(\hat y^1,\hat z^1)\|^q_{\s^q\times\M^q}.
  \end{equation*}
  We first assume that
  \[C_q\left(\intT{0}v^2(t)\dif t\right)^{q\over 2}<1.\]
  Since $\|(\hat y^1, \hat z^1)\|^q_{\s^q\times\M^q}<+\infty$, it follows immediately that
  $\{(y^n_t-y^1_t,z^n_t-z^1_t)_{t\in\tT}\}_{n=1}^{\infty}$ converges to some
  $(Y_t, Z_t)_{t\in\tT}$ in the space
  $\s^q\times \M^q$. Then we can deduce that $\{(y^n_t,z^n_t)_{t\in\tT}\}_{n=1}^{\infty}$
  converges to
  $(y_t:=Y_t+y^1_t,z_t:=Z_t+z^1_t)_{t\in\tT}$ in the space $\s^\beta\times \M^\beta$ for any
  $\beta\in(0,1)$ since $(y^1_t,z^1_t)_{t\in\tT}$ belongs to it.
  Moreover, since $(y^1_t)_{t\in\tT}$ belongs to the class (D) and the convergence
  in $\s^q$ with $q> 1$ is stronger than
  the convergence under the norm $\|\cdot\|_1$, we know that
  $\{(y^n_t)_{t\in\tT}\}_{n=1}^{\infty}$ converges to
  $(y_t)_{t\in\tT}$ under the norm $\|\cdot\|_1$. Then by taking the limit in
  both sides of BSDE \eqref{eq:BSDEYnPicardIterative} under ucp we can
  see that $(y_t,z_t)_{t\in\tT}$ is the desired solution of BSDE \eqref{eq:BSDEs}.

  For the general case, in view of $\intT{0}v^2(t)\dif t<+\infty$, we can as
  before subdivide the time interval $\tT$
  into a finite number of small intervals like $[T_i,T_{i+1}]$ such that
  \[
    C_q\left(\intT[T_{i+1}]{T_i}v^2(t)\dif t\right)^{q\over 2}<1.
  \]
  This completes the proof of the existence part of Theorem \ref{thm:MainResultL1Solution}.
\end{proof}

\begin{rmk}
  According to the proof procedure of the uniqueness and existence part of
  Theorem \ref{thm:MainResultL1Solution}, we know that if assumption
  \ref{H:gSublinearGrowthInZ} is satisfied as follows:
  \begin{equation*}
    |g(t,y,z)-g(t,y,0)|\leq \gamma(t)|z|^\alpha,
  \end{equation*}
  then $\gamma(t)$ in \ref{H:gSublinearGrowthInZ} need only to satisfy
  $\intT{0}\gamma^{2/(2-\alpha)}(t)\dif t<+\infty$.
\end{rmk}

\begin{rmk}\label{rmk:ComparisonWithBraindResults}
  In the case that $0\leq T<+\infty$, the $u(t)$, $v(t)$ and $\gamma(t)$ appearing
  in \ref{H:gMonotonicityInY}, \ref{H:gLipschitzInZ} and \ref{H:gSublinearGrowthInZ}
  are all bounded by a constant $c>0$ in the corresponding assumptions
  in \citet{BriandDelyonHu2003SPA}, and they do not deal with the case $T=+\infty$.
  But in our framework, $u(t)$, $v(t)$ and $\gamma(t)$ may be unbounded,
  and their integrability is the only requirement.
\end{rmk}

\begin{ex}\label{ex:L1UnboundedParameters}
  Let $0\leq T<+\infty$, $k=1$ and
  $$g(t,y,z)=\frac{1}{\sqrt[3]{t\,}}
  \big(\me^{-y}\one{y\leq 0}+(1-y^2)\one{y>0}\big)
  +\frac{t+1}{\sqrt[4]{t\,}}\big(|z|^2\wedge\sqrt{|z|}\,\big)+\frac{1}{1+t^4}.$$
  It is not hard to check that the generator $g$ satisfies assumptions
  \ref{H:OnlyIntegrable}, \ref{H:gContinuousInY} -- \ref{H:gSublinearGrowthInZ}
  with $\alpha=1/2$. Then Theorem \ref{thm:MainResultLpSolution} leads to that
  for each $\xi\in L^1(\Omega,\F_T,\PR;\rtn)$, BSDE $(\xi,T,g)$
  has a unique solution $(y_t,z_t)_{t\in\tT}$ in $\s^\beta\times\M^\beta$ for
  $\beta\in(1/2,1)$ with $(y_t)_{t\in\tT}$ belonging to the class (D).
  Remark \ref{rmk:ComparisonWithBraindResults} applies.
\end{ex}

\begin{ex}
  Let $0\leq T\leq +\infty$, $k=2$ and
  \begin{equation*}
    g(t,y,z)=\frac{1}{1+t^2}
    \begin{bmatrix}
      \displaystyle\me^{-y_1}+3y_2\\
      \displaystyle-\me^{y_2}-3y_1
    \end{bmatrix}+\me^{-t}
    \begin{bmatrix}
      \displaystyle\sin|z_1|\\[3pt]
      \displaystyle\sin|z_2|
    \end{bmatrix}+
    \begin{bmatrix}
      \displaystyle\me^{-t}\sin t\\[3pt]
      \displaystyle t\me^{-t}
    \end{bmatrix}.
  \end{equation*}
  where $y_i$ and $z_i$ $(i=1,2)$ stand for the $i$th component of the vector $y$
  and the $i$th row of the matrix $z$ respectively.
  It is not hard to verify that $g$ satisfies assumptions
  \ref{H:OnlyIntegrable}, \ref{H:gContinuousInY} -- \ref{H:gSublinearGrowthInZ}
  with $u(t)=1/(1+t^2)$, $v(t)=\gamma(t)=\me^{-t}$
  and for any $\alpha\in(0,1)$. Thus, by Theorem \ref{thm:MainResultLpSolution}
  we know that for each $\xi\in L^1(\Omega,\F_T,\PR;\rtn^k)$,
  BSDE $(\xi,T,g)$ has a unique solution $(y_t,z_t)_{t\in\tT}$ in
  $\s^\beta\times\M^\beta$ for $\beta\in(\alpha,1)$ with
  $(y_t)_{t\in\tT}$ belonging to the class (D).
  Remark \ref{rmk:ComparisonWithBraindResults} applies.
\end{ex}

\appendix
\section{Complement proofs of some lemmas}\label{sec:Appendix}
This appendix gives the detailed proofs of some lemmas.
\begin{proof}[\bf Proof of Lemma \ref{lem:estimate-z}]
  Let us fix the nonnegative function $\beta(t):\tT\mapsto\rtn^+$ with
  $\intT{0}\beta(t)\dif t<+\infty$ and
  $\beta(t)\geq 2\big(\mu(t)+\lambda^2(t)\big)$ for
  each $t\in\tT$. Similar to the change of variables in Remark \ref{rmk:ChangeOfVariables},
  we define $\overline{y}_t=\me^{\frac{1}{2}\intT[t]{0}\beta(s)\dif s}y_t$,
  $\overline{z}_t=\me^{\frac{1}{2}\intT[t]{0}\beta(s)\dif s}z_t$. Then $(\overline{y}_t,\overline{z}_t)_{t\in\tT}$
  solves BSDE $(\me^{\frac{1}{2}\intT{0}\beta(s)\dif s}\xi,T,\overline{g})$
  where
  $$\overline{g}(t,y,z):=\me^{\frac{1}{2}\intT[t]{0}\beta(s)\dif s}
    g\big(t,\me^{-\frac{1}{2}\intT[t]{0}\beta(s)\dif s}y,\me^{-\frac{1}{2}\intT[t]{0}\beta(s)\dif s}z\big)
    -\frac{1}{2}\beta(t)y,$$
  which satisfies assumption (A) with
  $\overline{\mu}(t)=\mu(t)-\frac{1}{2}\beta(t)$,
  $\overline{\lambda}(t)=\lambda(t)$ and
  $\overline{f}_t=\me^{\frac{1}{2}\intT[t]{0}\beta(s)\dif s}f_t$.
  The integrability conditions are equivalent with or without the superscript
  ``$\overline{{\color{white} a}}$" since $\intT{0}\beta(t)\dif t<+\infty$.
  Thus, with this change
  of variables we reduce to the case $\beta(t)\equiv 0$ and
  $\mu(t)+\lambda^2(t)\leq 0$. With omitting the superscript ``$\overline{{\color{white} a}}$"
  for notational convenience, we need to prove that there exists a constant
  $C^1_p>0$ such that for each $0\leq r\leq t\leq T$,
  \begin{equation}\label{eq:PriorEstimateOfZInAppendix}
    \EX\left[\left.
       \left(
         \intT{t}|z_s|^2\dif s
       \right)^{p\over 2}\right|\F_r
       \right]
    \leq C^1_p\EX
       \left[\left.
         \sup_{s\in\tT[t]}|y_s|^p
       +\left(\intT{t}f_s\dif s\right)^p
       \right|\F_r
       \right].
  \end{equation}

  For each $n\in\N$, let us introduce the following stopping time:
  $$\tau_n=\inf\left\{t\in\tT:\intT[t]{0}|z_s|^2\dif s\geq n\right\}\wedge T.$$
  Applying \Ito{} formula to
  $|y_s|^2$ yields that
  \begin{align}
    |y_{t\wedge\tau_n}|^2
    +\intT[\tau_n]{t\wedge\tau_n}|z_s|^2\dif s
    =|y_{\tau_n}|^2
    +2\intT[\tau_n]{t\wedge\tau_n}\langle y_s,g(s,y_s,z_s)\rangle\dif s
      -2\intT[\tau_n]{t\wedge\tau_n}\langle y_s,z_s\dif B_s\rangle.\label{eq:ito-result}
  \end{align}
  The inner product term including $g$ can be enlarged via assumption (A) (stated
  between Lemma \ref{lem:Tanaka-Briand} and Lemma \ref{lem:estimate-z}),
  $2ab\leq 2a^2+b^2/2$, $2ab\leq a^2+b^2$ and $\mu(t)+\lambda^2(t)\leq 0$
  as follows:
  \begin{align*}
    2\intT[\tau_n]{t\wedge\tau_n}\langle y_s,g(s,y_s,z_s)\rangle\dif s
    & \leq \intT[\tau_n]{t\wedge\tau_n}\left(2\mu(s)|y_s|^2+2\lambda^2(s)|y_s|^2+\frac{1}{2}|z_s|^2+2f_s|y_s|\right)\dif s\\
    & \leq \sup_{s\in [t\wedge\tau_n,\tau_n]}|y_s|^2
           +\frac{1}{2}\intT[\tau_n]{t\wedge\tau_n}|z_s|^2\dif s
           +\left(
              \intT[\tau_n]{t\wedge\tau_n}f_s\dif s
            \right)^2.
  \end{align*}
  Putting the previous inequality into \eqref{eq:ito-result} and noticing that
  $\tau_n\leq T$, we can deduce
  that there exists a constant $c_p>0$ depending only on $p$ such that for each
  $n\in\N$,
  \begin{align}
    \left(
      \intT[\tau_n]{t\wedge\tau_n}|z_s|^2\dif s
    \right)^{p\over 2}
    &\leq c_p
     \left[
        \sup_{s\in\tT[t\wedge\tau_n]}|y_s|^p+
       \left(
         \intT{t\wedge\tau_n}f_s\dif s
       \right)^p
     +\left|
        \intT[\tau_n]{t\wedge\tau_n}
        \langle y_s,z_s\dif B_s\rangle
      \right|^{p\over 2}
     \right].\label{eq:estimate-Z-before-expectation}
  \end{align}
  Moreover, the BDG inequality yields that there exists a constant $d_p>0$ depending
  only on $p$ such that for each $n\in\N$
  and $0\leq r\leq t\leq T$,
  \begin{align*}
    &c_p\EX
     \left[\left.\left|
       \intT[\tau_n]{t\wedge\tau_n}\langle y_s,z_s\dif B_s\rangle
       \right|^{p\over 2}\right|\F_r
     \right]
    \leq d_p\EX
             \left[\left.
               \left(
                 \intT[\tau_n]{t\wedge\tau_n}|y_s|^2|z_s|^2\dif s
               \right)^{p\over 4}\right|\F_r
             \right]\\
    &\leq \frac{d_p^2}{2}
          \EX
          \left[\left.
            \sup_{s\in\tT[t\wedge\tau_n]}|y_s|^p\right|\F_r
          \right]
          +\frac{1}{2}\EX
          \left[\left.
          \left(
            \intT[\tau_n]{t\wedge\tau_n}|z_s|^2\dif s
          \right)^{p\over 2}\right|\F_r
          \right].
  \end{align*}
  Thus, by taking conditional expectation with respect to $\F_r$ in
  both sides of the inequality
  \eqref{eq:estimate-Z-before-expectation} and then making use of
  Fatou's lemma,
  we can deduce that there exists a constant $C^1_p>0$ depending only on $p$ which
  satisfies estimate \eqref{eq:PriorEstimateOfZInAppendix}. This completes the
  proof of Lemma \ref{lem:estimate-z}.
\end{proof}

\begin{proof}[\bf Proof of Lemma \ref{lem:estimate-y}]
  Let the assumptions of Lemma \ref{lem:estimate-z} hold and $p>1$.
  Fix the nonnegative function $\beta(t):\tT\mapsto\rtn^+$ with
  $\intT{0}\beta(t)\dif t<+\infty$ and
  $\beta(t)\geq p\{\mu(t)+\lambda^2(t)/[1\wedge(p-1)]\}$. As in the proof
  of Lemma \ref{lem:estimate-z}, we also make the change of variables
  $\overline{y}_t=\me^{\frac{1}{2}\intT[t]{0}\beta(s)\dif s}y_t$,
  $\overline{z}_t=\me^{\frac{1}{2}\intT[t]{0}\beta(s)\dif s}z_t$.
  This reduces to the case $\beta(t)\equiv 0$ and
  $\mu(t)+\lambda^2(t)/[1\wedge(p-1)]\leq 0$. With omitting the superscript
  ``$\overline{{\color{white} a}}$", we have to prove that there exists a
  constant $C^2_p>0$ depending only on $p$ such that for each $0\leq r\leq t\leq T$,
  \begin{equation}\label{eq:PriorEstimateOfYInAppendix}
   \EX\left[\left.
       \sup_{s\in\tT[t]}|y_s|^p\right|\F_r
     \right]
   \leq C^2_p\EX
     \left[
       \left.|\xi|^p
       +\left(\intT{t}f_s\dif s\right)^p\right|\F_r
     \right].
  \end{equation}

  It follows from Lemma \ref{lem:Tanaka-Briand} that for each $t\in\tT$,
  \begin{align}\label{eq:tanaka-multi-exp}
    |y_t|^p+c(p)\intT{t}|y_s|^{p-2}\one{|y_s|\neq 0}|z_s|^2\dif s
    \leq&\ |\xi|^p+p\intT{t}|y_s|^{p-2}\one{|y_s|\neq 0}
      \langle y_s,g(s,y_s,z_s)\rangle\dif s\nonumber\\
    &-p\intT{t}|y_s|^{p-2}\one{|y_s|\neq 0}\langle y_s,z_s\dif B_s\rangle.
  \end{align}
  In view of assumption (A), we can get that
  \begin{align*}
    c(p)\intT{0}|y_s|^{p-2}\one{|y_s|\neq 0}|z_s|^2\dif s
    \leq
    &\ |\xi|^p+p\intT{0}
       \big(
         \mu(s)|y_s|^p+\lambda(s)|y_s|^{p-1}|z_s|+f_s|y_s|^{p-1}
       \big)\dif s\\
    &-p\intT{0}|y_s|^{p-2}\one{|y_s|\neq 0}\langle y_s,z_s\dif B_s\rangle,
  \end{align*}
  from which we have $\prs$,
  \begin{equation*}
    \intT{0}|y_s|^{p-2}\one{|y_s|\neq 0}|z_s|^2\dif s<+\infty.
  \end{equation*}

  Now enlarge the inner product term
  including $g$ in \eqref{eq:tanaka-multi-exp} with $ab\leq \alpha a^2/4+b^2/\alpha$ $(\alpha=(p-1)\wedge 1)$
  and $\mu(t)+\lambda^2(t)/[1\wedge(p-1)]\leq 0$ as follows:
  \begin{align*}
    &p\intT{t}|y_s|^{p-2}\one{|y_s|\neq 0}\langle y_s,g(s,y_s,z_s)\rangle\dif s\\
    &\leq\intT{t}
     \left[
       p\mu(s)|y_s|^p+\frac{p\lambda^2(s)}{1\wedge(p-1)}|y_s|^p
       +\frac{c(p)}{2}|y_s|^{p-2}\one{|y_s|\neq 0}|z_s|^2
       +p|y_s|^{p-1}f_s
     \right]\dif s\\
    &\leq \intT{t}
     \left[
       \frac{c(p)}{2}|y_s|^{p-2}\one{|y_s|\neq 0}|z_s|^2
       +p|y_s|^{p-1}f_s
     \right]\dif s.
  \end{align*}
  Putting the previous inequality into \eqref{eq:tanaka-multi-exp} we can get
  that for each $t\in\tT$,
  \begin{equation}\label{eq:estimate-Y-withX}
    |y_t|^p
    +\frac{c(p)}{2}\intT{t}|y_s|^{p-2}\one{|y_s|\neq 0}|z_s|^2\dif s
    \leq X_t-p\intT{t}|y_s|^{p-2}\one{|y_s|\neq 0}\langle y_s,z_s\dif B_s\rangle,
  \end{equation}
  where $X_t:=|\xi|^p+p\intT{t}|y_s|^{p-1}f_s\dif s$.
  Note that $\{M_t:=\intT[t]{0}|y_s|^{p-2}\one{|y_s|\neq 0}\langle y_s,z_s\dif B_s\rangle\}_{t\in\tT}$ is
  a uniformly integrable martingale. In fact, it follows from the BDG inequality
  and Young's inequality (for any
  nonnegative constant $a$ and $b$, $ab\leq a^p/p+b^q/q$ holds true with
  $q=p/(p-1)$) that
  \begin{align*}
    \EX
    \left[
      \langle M,M\rangle^{1/2}_T
    \right]
    & \leq
    \EX
    \left[
      \left(
        \intT{0}|y_s|^{2p-2}|z_s|^2\dif s
      \right)^{1\over 2}
    \right]
    \leq \EX
      \left[
        \sup_{s\in\tT}|y_s|^{p-1}\cdot
        \left(
          \intT{0}|z_s|^2\dif s
        \right)^{1\over 2}
      \right]\\
    & \leq \frac{p-1}{p}\EX
      \left[
        \sup_{s\in\tT}|y_s|^p
      \right]+{1\over p}\EX
      \left[
        \left(
          \intT{0}|z_s|^2\dif s
        \right)^{p\over 2}
      \right].
  \end{align*}
  Since $(y_t)_{t\in\tT}$ belongs to $\s^p$ and then $(z_t)_{t\in\tT}$ belongs to
  $\M^p$ by Lemma \ref{lem:estimate-z}, the right term in the previous inequality
  is finite.
  Therefore, for each $0\leq r\leq t\leq T$
  the inequality \eqref{eq:estimate-Y-withX} yields both
  \begin{equation*}
    \frac{c(p)}{2}\EX
    \left[\left.
      \intT{t}|y_s|^{p-2}\one{|y_s|\neq 0}|z_s|^2\dif s\right|\F_r
    \right]\leq\EX[X_t|\F_r],
  \end{equation*}
  and
  \begin{align*}
    \EX
    \left[\left.
      \sup_{s\in\tT[t]}|y_s|^p\right|\F_r
    \right]
   &\leq \EX[X_t|\F_r]+
    p\EX\left[\left.\sup_{s\in\tT[t]}
        \left|
          \intT[T]{s}|y_{\bar u}|^{p-2}\one{|y_{\bar u}|\neq 0}\langle y_{\bar u},z_{\bar u}\dif B_{\bar u}\rangle
        \right|\right|\F_r\right]\\
   &\hspace{-2cm}\leq \EX[X_t|\F_r]+k_p\EX
    \left[\left.
      \left(\intT{t}|y_s|^{2p-2}\one{|y_s|\neq 0}|z_s|^2\dif s\right)^{1\over 2}\right|\F_r
    \right]\\
   &\hspace{-2cm}\leq \EX[X_t|\F_r]+\frac{1}{2}\EX
    \left[\left.
      \sup_{s\in\tT[t]}|y_s|^p\right|\F_r
    \right]
    +\frac{k_p^2}{2}\EX
     \left[\left.
       \intT{t}|y_s|^{p-2}\one{|y_s|\neq 0}|z_s|^2\dif s\right|\F_r
     \right],
  \end{align*}
  where the constant $k_p>0$ depending only on $p$ follows by the BDG inequality.
  Combining the previous two inequalities
  we can get that for each $t\in\tT$,
  \begin{equation}\label{eq:estimate-Y-with-EX}
    \EX
    \left[\left.
      \sup_{t\in\tT[t]}|y_s|^p\right|\F_r
    \right]
    \leq 2\left[1+\frac{k^2_p}{c(p)}\right]\EX[X_t|\F_r].
  \end{equation}
  Let $l_p:=2p\big(1+k^2_p/c(p)\big)$. It follows from Young's inequality that
  \begin{align*}
     l_p\EX
     \left[\left.
       \intT{t}|y_s|^{p-1}f_s\dif s\right|\F_r
     \right]
    &\leq \EX
     \left[\left.
       \left(\frac{p}{2p-2}\right)^{p-1\over p}
       \sup_{s\in\tT[t]}|y_s|^{p-1}\cdot
       l_p\left(\frac{p}{2p-2}\right)^{1-p\over p}\intT{t}f_s\dif s
       \right|\F_r
     \right]\\
    &\leq \frac{1}{2}\EX
     \left[\left.
       \sup_{s\in\tT[t]}|y_s|^p\right|\F_r
     \right]
     +K_p\EX
      \left[\left.
        \left(
          \intT{t}f_s\dif s
        \right)^p\right|\F_r
      \right],
  \end{align*}
  where $K_p:=(l_p)^p\big(p/(2p-2)\big)^{1-p}/p$. Now combining
  \eqref{eq:estimate-Y-with-EX} with the definition of $X_t$ and the previous
  inequality, we deduce that there must exist a constant $C^2_p>0$ depending
  only on $p$ such that \eqref{eq:PriorEstimateOfYInAppendix} holds true.
  The proof of Lemma \ref{lem:estimate-y} is then completed.
\end{proof}

\begin{proof}[\bf Proof of Lemma \ref{lem:AppendixProofOfgnWithLocallyLip}]
  It is clear that $f_n(t,y)$ satisfies assumption \ref{H:gContinuousInY}.
  By assumption \ref{H:gGeneralGrowthInY} on $g$, \eqref{eq:g-n-convolution}
  and \eqref{eq:g-n-convolution-equivalent}
  we can obtain that
  \begin{align*}
    |f_n(t,y)|&\leq \int_{\rtn^k}\rho_n(x)|f(t,y-x)|\dif x
               \leq\int_{\rtn^k}\rho_n(x)\big(|f(t,0)|+u(t)\varphi(|y-x|)\big)\dif x\\
              &\leq |f(t,0)|+u(t)\int_{\{x:|x|\leq 1\}}\rho(x)\varphi(|y-\frac{x}{n}|)\dif x
               \leq |f(t,0)|+u(t)\varphi(|y|+1).
  \end{align*}
  Thus, $f_n(t,y)$ satisfies \ref{H:gGeneralGrowthInY} with $\varphi$
  replaced by $\phi$, and then \eqref{eq:FnEnlargedByEAndU} follows from
  \eqref{eq:AssumptionFirstStep}.

  Furthermore, for each $y_1,y_2\in\rtn^k$, we have, in view of
  \ref{H:gMonotonicityInYWithU=0} on $g$,
  \begin{equation*}
    \langle y_1-y_2,f_n(t,y_1)-f_n(t,y_2)\rangle
    =\int_{\rtn^k}\rho_n(x)\langle y_1-y_2, f(t,y_1-x)-f(t,y_2-x)\rangle\dif x\leq 0,
  \end{equation*}
  which means that \ref{H:gMonotonicityInYWithU=0} holds true for $f_n$.

  Finally, fix $(\omega,t)\in\Omega\times\tT$, for the gradient of $f_n(t,y)$
  with respect to $y$, we have
  \begin{equation*}
    |\nabla f_n(t,y)|
    \leq\int_{\rtn^k}|\nabla\rho_n(y-x)||f(t,x)|\dif x,\quad \forall y\in\rtn^k.
  \end{equation*}
  It then follows from assumption \ref{H:gGeneralGrowthInY} on $g$ and
  \eqref{eq:AssumptionFirstStep} that for each $y\in\rtn^k$,
  \begin{align*}
    |\nabla f_n(t,y)|
    &\leq \int_{\rtn^k}|\nabla\rho_n(y-x)|
       \big(|f(t,0)|+u(t)\varphi(|x|)\big)\dif x\\
    &=\int_{\{x:|y-x|\leq 1/n\}}|\nabla\rho_n(y-x)|
       \big(|f(t,0)|+u(t)\varphi(|x|)\big)\dif x\\
    &\leq \big(K\me^{-t}+u(t)\varphi(|y|+1)\big)
       \int_{\rtn^k}|\nabla\rho_n(x)|\dif x.
  \end{align*}
  Then \eqref{eq:fnLocallyLipschitzContinuousInY} follows immediately.
  That is, $f_n(t,y)$ is locally Lipschitz continuous in
  $y$ non-uniformly with respect to $t$. Lemma \ref{lem:AppendixProofOfgnWithLocallyLip}
  is then proved.
\end{proof}

\begin{proof}[\bf Proof of Lemma \ref{lem:AppendixWeakConvergenceU=f}]
  For any real number $\varepsilon>0$, we set
  \begin{equation*}
    X^\varepsilon_t:=y_t-\frac{\varepsilon\big(U_t-f(t,y_t)\big)}{|U_t-f(t,y_t)|}\one{|U_t-f(t,y_t)|\neq 0}.
  \end{equation*}
  In view of \eqref{eq:YAndEXYInfinite}, it is clear that $\pts$,
  $|X^\varepsilon_t|\leq |y_t|+\varepsilon\leq a+\varepsilon$. It then follows from
  \eqref{eq:g-n-convolution-equivalent} that $\pts$,
  $f_n(t,X^\varepsilon_t)\to f(t,X^\varepsilon_t)$ as $n\to\infty$
  and from \eqref{eq:AssumptionFirstStep} and assumption
  \ref{H:gGeneralGrowthInY} on $f_n$ and $f$ that
  \begin{equation*}
    |f_n(t,X^\varepsilon_t)-f(t,X^\varepsilon_t)|
    \leq 2K\me^{-t}+u(t)\big(\phi(0)+\varphi(a+\varepsilon)+\phi(a+\varepsilon)\big).
  \end{equation*}
  Then Lebesgue's dominated convergence theorem leads to
  \begin{equation}\label{eq:lim-fn-f-Xt}
    \lim_{n\to\infty}\EX\left[\intT{0}|f_n(t,X^\varepsilon_t)-f(t,X^\varepsilon_t)|\dif t\right]= 0.
  \end{equation}
  Furthermore, we can deduce that
  \begin{equation}\label{eq:limsup-yn-Xt}
    \limsup_{n\to\infty}\EX
    \left[
      \intT{0}\langle y^n_t-X^\varepsilon_t,f_n(t,y^n_t)-f(t,X^\varepsilon_t)\rangle\dif t
    \right]\leq 0.
  \end{equation}
  Indeed, for each $n\in\N$, by \ref{H:gMonotonicityInYWithU=0} on $f_n$ we know
  that $\pts$,
  \begin{align*}
    \langle y^n_t-X^\varepsilon_t,f_n(t,y^n_t)-f(t,X^\varepsilon_t)\rangle\leq
    \langle y^n_t-X^\varepsilon_t,f_n(t,X^\varepsilon_t)-f(t,X^\varepsilon_t)\rangle.
  \end{align*}
  Then by the previous inequality and \eqref{eq:lim-fn-f-Xt} and noticing that
  $\pts$, $|y^n_t-X^\varepsilon_t|\leq 2a+\varepsilon$, we have
  \begin{align*}
   &\limsup_{n\to\infty}\EX
   \left[
     \intT{0}\langle y^n_t-X^\varepsilon_t,f_n(t,y^n_t)-f(t,X^\varepsilon_t)\rangle\dif t
   \right]\\
   &\leq(2a+\varepsilon)\limsup_{n\to\infty}\EX
    \left[
      \intT{0}|f_n(t,X^\varepsilon_t)-f(t,X^\varepsilon_t)|\dif t
    \right]=0.
  \end{align*}

  In the sequel, applying \Ito{} formula to $|y^n_t|^2$ we can get that
  \begin{equation*}
    2\EX
    \left[
      \intT{0}\langle y^n_t,f_n(t,y^n_t)\rangle\dif t
    \right]=|y_0^n|^2-\EX\left[|\xi|^2\right]
    +\EX\left[\intT{0}|z^n_t|^2\dif t\right].
  \end{equation*}
  Then since the mapping $z\mapsto\EX[\intT{0}|z_t|^2\dif t]$ is weakly lower
  semi-continuous and $y^n_0$ converges to $y_0$ in $\rtn^k$, we have
  \begin{equation}
    \liminf_{n\to\infty}2\EX
    \left[
      \intT{0}\langle y^n_t,f_n(t,y^n_t)\rangle\dif t
    \right]
    \geq |y_0|^2-\EX\left[|\xi|^2\right]+\EX\left[\intT{0}|z_t|^2\dif t\right]
    =2\EX
     \left[
       \intT{0}\langle y_t,U_t\rangle\dif t
     \right].\label{eq:LiminfYnWeakConvergenceWithU}
  \end{equation}
  The equal sign in the previous equation follows from applying
  \Ito{} formula to $|y_t|^2$.
  Combining the weak convergences with \eqref{eq:LiminfYnWeakConvergenceWithU}
  and \eqref{eq:limsup-yn-Xt}, we can deduce that
  \begin{equation*}
    \EX
    \left[
      \intT{0}\langle y_t-X^\varepsilon_t, U_t-f(t,X^\varepsilon_t)\rangle\dif t
    \right]
    \leq \liminf_{n\to\infty}\EX
    \left[
      \intT{0}\langle y^n_t-X^\varepsilon_t,f_n(t,y^n_t)-f(t,X^\varepsilon_t)\rangle\dif t
    \right]\leq 0.
  \end{equation*}
  Thus, noticing the definition of $X^\varepsilon_t$, we have, for each $\varepsilon>0$,
  \begin{equation*}
    \EX
    \left[
      \intT{0}\left\langle \frac{U_t-f(t,y_t)}{|U_t-f(t,y_t)|}\one{|U_t-f(t,y_t)|\neq 0},
      U_t-f(t,X^\varepsilon_t)\right\rangle \dif t
    \right]\leq 0.
  \end{equation*}
  Sending $\varepsilon$ to $0$ yields that $\pts$, $X^\varepsilon_t\to y_t$. Then
  noticing that \ref{H:gContinuousInY} holds true for $f$, we have $\pts$,
  $f(t,X^\varepsilon_t)\to f(t,y_t)$. Moreover, since $\EX[(\intT{0}|U_t|\dif t)^2]<+\infty$
  and $|f(t,X^\varepsilon_t)|\leq K\me^{-t}+u(t)\varphi(a+\varepsilon)$,
  Lebesgue's dominated convergence theorem leads to that
  \begin{equation*}
    \EX
    \left[
      \intT{0}|U_t-f(t,y_t)|\dif t
    \right]\leq 0,
  \end{equation*}
  from which we have $\pts$, $U_t=f(t,y_t)$.
\end{proof}

\begin{proof}[\bf Proof of Lemma \ref{lem:AppendixHnProperty}]
  It follows from the definition of $h_n$ that $h_n$ satisfies
  \ref{H:gContinuousInY} and $h_n(t,0,V_t)=g(t,0,V_t)$. Therefore,
  $|h_n(t,0,V_t)|\leq K\me^{-t}$.
  Next we check that $h_n$ satisfies \ref{H:gMonotonicityInY}. For each
  $y_1$, $y_2\in\rtn^k$, if $|y_1|>r'+1$ and $|y_2|>r'+1$, \ref{H:gMonotonicityInY}
  is trivially satisfied and thus we reduce to the case where $|y_2|\leq r'+1$.
  For notation convenience, we set, for each $n\in\N$ and $t\in\tT$,
  \begin{equation}\label{eq:PinPsin}
    \pi^{n}(t):=\pi_{n\me^{-t}}(V_t),\quad
    \psi^{n}(t):=\frac{n\me^{-t}}{\psi_{r'+1}(t)\vee(n\me^{-t})}.
  \end{equation}
  By adding and subtracting $\theta_{r'}(y_1)g(t,y_2,\pi^n(t))$
  we can deduce that
  \begin{align*}
    \langle y_1-y_2, h_n(t,y_1,V_t)-h_n(t,y_2,V_t)\rangle
    & =\psi^n(t)\theta_{r'}(y_1)
       \langle y_1-y_2,g(t,y_1,\pi^n(t))-g(t,y_2,\pi^n(t))\rangle\\
    & \quad+\psi^n(t)
      \big(\theta_{r'}(y_1)-\theta_{r'}(y_2)\big)
      \langle y_1-y_2,g(t,y_2,\pi^n(t))-g(t,0,\pi^n(t))\rangle.
  \end{align*}
  The first term on the right side is non-positive since $g$ satisfies
  \ref{H:gMonotonicityInYWithU=0}. For the second term, since
  $\theta_{r'}$ is $C(r')$-Lipschitz and $|y_2|\leq r'+1$, we can get that
  \begin{align*}
    \big(\theta_{r'}(y_1)-\theta_{r'}(y_2)\big)
     \langle y_1-y_2,g(t,y_2,\pi^n(t))-g(t,0,\pi^n(t))\rangle
    \leq C(r')|y_1-y_2|^2|g(t,y_2,\pi^n(t))-g(t,0,\pi^n(t))|.
  \end{align*}
  And assumption \ref{H:gLipschitzInZ} on $g$ and the definition
  of $\psi_{r'}(t)$ in \ref{H:gGeneralizedGeneralGrowthInY} yield that
  \begin{align}
    & |g(t,y_2,\pi^n(t))-g(t,0,\pi^n(t))|
      \leq v(t)\pi^n(t)+|g(t,y_2,0)-g(t,0,\pi^n(t))|\nonumber\\
    & \leq 2v(t)\pi^n(t)+|g(t,y_2,0)-g(t,0,0)|
      \leq 2v(t)n\me^{-t}+\psi_{r'+1}(t).\label{eq:gEnlargeWithPsiHnProperty}
  \end{align}
  Then, in view of $|\psi^n(t)|\leq 1$ and
  $|\psi^n(t)\psi_{r'+1}(t)|=|\pi_{n\me^{-t}}(\psi_{r'+1}(t))|\leq n\me^{-t}$, we have
  \begin{align*}
    \langle y_1-y_2, h_n(t,y_1,V_t)-h_n(t,y_2,V_t)\rangle
    \leq nC(r')\big(2v(t)\me^{-t}+\me^{-t}\big)|y_1-y_2|^2.
  \end{align*}
  Thus, note that $2\intT{0}v(t)\me^{-t}\dif t\leq \intT{0}v^2(t)\dif t+\intT{0}\me^{-2t}\dif t<+\infty$,
  we know that $h_n$ satisfies \ref{H:gMonotonicityInY}.

  Finally, we check that $h_n$ satisfies \ref{H:gGeneralGrowthInY}.
  It follows from the definition of $h_n$ that
  \begin{align*}
    |h_n(t,y,V_t)|
    \leq \psi^n(t)\theta_{r'}(y)
         |g(t,y,\pi^n(t))-g(t,0,\pi^n(t))|+|h_n(t,0,V_t)|.
  \end{align*}
  And similar to \eqref{eq:gEnlargeWithPsiHnProperty}, we know that for each
  $y\in\rtn^k$ with $|y|\leq r'+1$,
  \begin{equation*}
    |g(t,y,\pi^n(t))-g(t,0,\pi^n(t))|
    \leq 2nv(t)\me^{-t}+\psi_{r'+1}(t).
  \end{equation*}
  Hence, we deduce that for each $y\in\rtn^k$,
  \begin{equation*}
    |h_n(t,y,V_t)|\leq |h_n(t,0,V_t)|+ n\big(2v(t)\me^{-t}+\me^{-t}\big),
  \end{equation*}
  which means that \ref{H:gGeneralGrowthInY} is satisfied for $h_n$.
\end{proof}

\begin{proof}[\bf Proof of Lemma \ref{lem:AppendixHnEnlargeForLebesgueConvrgence}]
  We will use the same notations in \eqref{eq:PinPsin}.
Let
  \begin{align*}\allowdisplaybreaks
    A(s)&=\big(
            g(s,y^n_s,\pi^{n+i}(s))-g(s,y^n_s,\pi^{n}(s))
          \big)\psi^{n+i}(s),\\
    B(s)&=\big(
            g(s,y^n_s,\pi^{n}(s))-g(s,y^n_s,0)
          \big)
          \big(
            \psi^{n+i}(s)-\psi^n(s)
          \big),\\
    C(s)&=\big(
            g(s,y^n_s,0)-g(s,0,0)
          \big)
          \big(
            \psi^{n+i}(s)-\psi^n(s)
          \big),\\
    D(s)&=\big(
            g(s,0,0)-g(s,0,\pi^n(s))
          \big)
          \big(
            \psi^{n+i}(s)-\psi^n(s)
          \big),\\
    E(s)&=\big(
            g(s,0,\pi^n(s))-g(s,0,\pi^{n+i}(s))
          \big)\psi^{n+i}(s).
  \end{align*}
  Then we have
  \[
    |h'_{n+i}(s,y^n_s,V_s)-h'_n(s,y^n_s,V_s)|=|A(s)+B(s)+C(s)+D(s)+E(s)|.
  \]
  It follows from assumption \ref{H:gLipschitzInZ}
  on $g$, the fact $|\psi^{n+i}(s)|\leq 1$ and the definitions of $\pi^n$ and
  $\pi^{n+i}$ that
  \begin{equation*}
    |A(s)|
    \leq v(s)|\pi^{n+i}(s)-\pi^n(s)|
      =v(s)|V_s|
      \left|
        \frac{(n+i)\me^{-s}}{|V_s|\vee\big((n+i)\me^{-s}\big)}-
        \frac{n\me^{-s}}{|V_s|\vee (n\me^{-s})}
      \right|
    \leq v(s)|V_s|\one{|V_s|>n\me^{-s}}.
  \end{equation*}
  Similar to the previous proof procedure, we also have that
  $|E(s)|\leq v(s)|V_s|\one{|V_s|>n\me^{-s}}$.

  Next, in view of assumption \ref{H:gLipschitzInZ} on $g$
  and the fact $|\psi^{n+i}(s)-\psi^n(s)|\leq \one{\psi_{r'+1}(s)>n\me^{-s}}$,
  we can get
  \begin{equation*}
    |B(s)|\leq v(s)|\pi^n(s)|\one{\psi_{r'+1}(s)>n\me^{-s}}
    \leq v(s)|V_s|\one{\psi_{r'+1}(s)>n\me^{-s}}.
  \end{equation*}
  Similarly, we also have $|D(s)|\leq v(s)|V_s|\one{\psi_{r'+1}(s)>n\me^{-s}}$.

  Finally, since $\pts$, $|y^n_s|\leq r'$, it follows from
  assumption \ref{H:gGeneralizedGeneralGrowthInY} that
  \begin{equation*}
    |C(s)|\leq |g(s,y^n_s,0)-g(s,0,0)||\psi^{n+i}(s)-\psi^n(s)|
    \leq \psi_{r'+1}(s)\one{\psi_{r'+1}(s)>n\me^{-s}}.
  \end{equation*}
  Then we can obtain the result as follows:
  \begin{equation*}
    |h'_{n+i}(s,y^n_s,V_s)-h'_n(s,y^n_s,V_s)|
    \leq 2v(s)|V_s|\one{|V_s|>n\me^{-s}}+2v(s)|V_s|\one{\psi_{r'+1}(s)>n\me^{-s}}
         +\psi_{r'+1}(s)\one{\psi_{r'+1}(s)>n\me^{-s}},
  \end{equation*}
  which completes this proof.
\end{proof}


\renewcommand{\baselinestretch}{1.3}
\setlength{\bibsep}{2pt}

\end{document}